\documentclass[12pt,journal,onecolumn ]{IEEEtran}

\IEEEoverridecommandlockouts                              

 \usepackage{datetime}  
\usepackage{hyperref}       
\hypersetup{colorlinks=true, linkcolor=blue, breaklinks=true, urlcolor=blue}

\IEEEoverridecommandlockouts              

\usepackage{amsmath,amssymb}
\usepackage{amsthm,thmtools}
\usepackage{amsfonts} 
\usepackage{blkarray}
 \usepackage{graphicx}
\usepackage{enumerate}
\usepackage{color}                  

\usepackage{verbatim}
\usepackage{bbm}

\usepackage{amssymb}
\usepackage{amsbsy}
\usepackage{amsmath,lipsum}
\usepackage{setspace}
\usepackage{url}

\usepackage{float}
 \usepackage{dsfont}
 \usepackage{cuted} \stripsep-3pt

 \usepackage{caption} 
 \newcommand{\Int}{\operatorname{int}}
 
 \newcommand{\real}{\operatorname{Re}}
 
 \newcommand{\diag}{\operatorname{diag}}
 
 \newcommand{\tr}{\operatorname{trace}}
 \newcommand{\adju}{\operatorname{adj}}
\newcommand{\st}{:}

\newcommand{\spanop}{\operatorname{span}}

\newcommand{\trace}{\operatorname{trace}}

\declaretheorem[name={Example},qed={\lower-0.3ex\hbox{$\square$}} ] {Example}

\declaretheorem[name={Definition}  ] {Definition}
\declaretheorem[name={Theorem}  ] {Theorem}
\declaretheorem[name={Lemma}  ] {Lemma}
\declaretheorem[name={Remark}  ] {Remark}
\declaretheorem[name={Corollary}  ] {Corollary}

\declaretheorem[name={Proposition}  ] {Proposition}

\newcommand {\R}{\mathbb R}

\newcommand {\C}{\mathbb C}

\newcommand{\be}{\begin{equation}}
\newcommand{\ee}{\end{equation}}
\newcommand{\sgn}{\operatorname{{\mathrm sgn}}}

\newcommand{\U}{\mathcal U}
\newcommand{\N}{\mathbb N_0}
\newcommand{\B}{\mathcal B \mathcal B}

\newcommand*\dif{\mathop{}\!\mathrm{d}}

\usepackage{lineno}

\DeclareMathOperator{\vol}{vol}
\DeclareMathOperator{\cof}{cof}

 \title{\LARGE \bf
Compound matrices in systems and control theory:\\ a tutorial \\  \vspace*{20pt} \normalsize  \today{}  }
 
\author{Eyal Bar-Shalom and Omri Dalin and   Michael Margaliot
\thanks{This research was  partially supported by   research grants from the~ISF and~DFG.}
\thanks{An abridged version of this paper has been presented
in 
  the 2021 IEEE Conf. on Decision and Control~\cite{comp_barshalom}.\\
The authors  are  with 
 the School of Electrical  Engineering,
		Tel-Aviv University, Tel-Aviv~69978, Israel.
		Corresponding author: Michael Margaliot (E-mail: \texttt{michaelm@tauex.tau.ac.il)
		}}}
\begin{document}
\maketitle
%
\begin{center}
Dedicated to Eduardo D. Sontag, a friend and mentor,   on the occasion of
  his 70th birthday  \vspace*{0.75cm}
\end{center}
\begin{abstract}
%
The  multiplicative and  additive compounds of a matrix play an important role in several fields of mathematics including geometry, 
multi-linear algebra, combinatorics, and the  analysis of
nonlinear time-varying dynamical  systems. There is a growing interest in applications of these compounds, and their generalizations, in systems and control theory. 
The goal of this  tutorial paper is to provide a gentle and self-contained  introduction to these topics with an emphasis on  the geometric interpretation of the compounds, and  to describe 
some of their recent applications including several  non-trivial 
generalizations of positive systems, cooperative systems, contracting systems, and more. 
%
\end{abstract}
%
\begin{IEEEkeywords}
 Contracting systems, diagonal stability, positive systems, cooperative systems, chaotic systems, sign variation diminishing property,  volume of  parallelotopes, Hankel $k$-positive systems.
\end{IEEEkeywords}
\newpage 

\tableofcontents
\newpage 
\section{Introduction}
%
 Let~$A\in\C^{n\times n}$. Fix~$k\in\{1,\dots,n\}$.
 The $k$-multiplicative and $k$-additive compounds of~$A$, denoted~$A^{(k)}$ and~$A^{[k]}$, respectively,  are~$\binom{n}{k}\times
 \binom{n}{k}$ matrices that 
 play an important role in several fields of applied mathematics.
These matrices have an interesting spectral property. If~$\lambda_i$, $i=1,\dots,n$, denote
the eigenvalues of~$A$
then the eigenvalues of~$A^{(k)}$ are the~$\binom{n}{k}$ products:
\be\label{eq:prod}
\bigg\{ \prod_{j=1}^k \lambda_{i_j}\st 1\leq i_1<\dots<i_k\leq n\bigg\},
\ee
and those of~$A^{[k]}$ are the~$\binom{n}{k}$ sums:
 \be\label{eq:sumeig}
\bigg\{\sum_{j=1}^k \lambda_{i_j}\st 1\leq i_1<\dots<i_k\leq n\bigg\}.
 \ee
 In particular,~$A^{(1)}=A^{[1]}=A$, so the eigenvalues of these matrices are just the~$\lambda_i$s, 
 $A^{(n)}=\det(A)=\prod_{i=1}^n \lambda_i$, and~$A^{[n]}=\trace(A)=\sum_{i=1}^n \lambda_i$.

 Recently, there is a growing interest in the applications of these compounds, and their generalizations in systems and control theory (see, e.g.~\cite{cheng_diag_stab,kordercont,rami_osci,rola_spect,CTPDS,grussler2021internally,gruss2,gruss3,grussler2022variation,gruss5, margaliot2019revisiting,wu2020generalization,Eyal_k_posi,DT_K_POSI,9107214,TPDS_NEAR_EQ}).
 This tutorial paper reviews the~$k$-compounds,  focusing on their  geometric interpretation, and
 surveys some of their  recent applications in systems and control theory,  including the introduction and analysis of~$k$-positive systems,
 $k$-cooperative systems, 
 $k$-contracting systems, 
    $k$-diagonal stability, and Hankel $k$-positive input/output systems.

This paper is organized as follows: The next section provides some geometric  motivation 
by relating the evolution of the volume of parallelotopes under a linear time invariant~(LTI) system and sums of eigenvalues such as in~\eqref{eq:sumeig}.  
Section~\ref{sec:compund_matrices}   introduces the multiplicative     compound and reviews some of its  properties and, in particular, 
  the fundamental role of the $k$-multiplicative compound 
in computing the volume of~$k$-dimensional parallelotopes,
and in establishing sign variation diminishing properties.
Section~\ref{sec:add_comp} describes the additive compound, and 
this  sets  the stage to explaining the role of the   compounds  in ordinary differential equations~(ODEs) in Section~\ref{sec:odes}. The following  section describes a general and important principle for what we call  
 $k$-generalizations of dynamical systems. The idea is to take a dynamical property, e.g. contraction in a nonlinear system, and require that it holds for $k$-dimensional bodies rather than~$1$-dimensional bodies (i.e., lines). For~$k=1$ this reduces to standard contraction. We then demonstrate  how this general principle leads to interesting  and non-trivial generalizations of contracting, positive, cooperative, and diagonally stable systems to
 $k$-contracting, $k$-positive, $k$-cooperative, and~$k$-diagonally stable systems. Section~\ref{sec:alpha}
 reviews the  recently introduced 
 concept of~$\alpha$-compounds, with~$\alpha$ being a \emph{real} number in~$[1,n]$. For~$\alpha \in (k,k+1)$, the $\alpha$-compound can be interpreted as an
 interpolation between the~$k$  and~$k+1$ compounds. We show that this leads to the notion of~$\alpha$-contracting systems. Any attractor of such a  system has a Hausdorff dimension smaller than~$\alpha$. 
 Section~\ref{sec:hankel} reviews  a type 
 of~$k$-positivity in systems with an input and output. 
  The final section concludes and describes several directions for future research. 

\subsection{Notation}
We use standard notation. The positive orthant in~$\R^n$ is~$\R^n_+ := \{x\in\R^n: x_i\geq 0,\; i=1,\dots,n\}$, and its interior is~$\R^n_{++} := \{x\in\R^n: x_i >  0,\; i=1,\dots,n\}$.
$\mathbb{Z}$ denotes the set of integers, and~$\N$ is the 
subset of non-negative integers. For a set~$S$, 
$\Int(S)$ denotes  the interior of~$S$.
For scalars~$\lambda_i$, $i\in\{1,\dots, n\}$, $\diag(\lambda_1,\dots,\lambda_n)$ is the $n\times n$ diagonal matrix with
diagonal entries~$\lambda_i$. 
Vectors [matrices] are denoted by small [capital] letters. For a matrix~$A$, $A^T$ is the transpose of~$A$. For a square matrix~$B$, $\tr(B)$ [$\det(B)$] is the trace [determinant]  of~$B$. Inequalities between  matrices~$A,B\in\R^{n\times m}$ are interpreted component-wise, e.g.~$A\geq B$ if~$a_{ij}\geq b_{ij}$ for every~$i,j$, and~$A\gg B$ if~$a_{ij}>b_{ij}$
for every~$i,j$.
A matrix is called Metzler if all its off-diagonal entries are non-negative. 
 A matrix is called Toeplitz [Hankel] if the entries along any diagonal [anti-diagonal] are equal. A matrix is called Hurwitz if the real part of every eigenvalue of the matrix is negative. A matrix is called Schur if the absolute value  of every eigenvalue of the matrix is smaller than~$1$. A      square binary matrix that has exactly one entry of~$1$ in each row and each column, and~$0$s elsewhere
is called a permutation matrix.  A matrix~$A\in\R^{n\times n} $ is called reducible if there  exists a permutation matrix~$P\in\{0,1\}^{n\times n} $ such that~$PAP^T=\begin{bmatrix}B&C\\0&D\end{bmatrix}$, where~$0$ is an~$r\times (n-r)$ zero block matrix for some~$1\leq r\leq n-1$. A matrix is called irreducible if it is not reducible. 
 
 Compound matrices require notation for specifying the minors of a matrix. 
 Let~$Q(k,n)$ denote all
 the~$\binom{n}{k}$
increasing sequences of~$k$ integers from the set~$\{1,\dots,n\}$, ordered lexicographically.
For example,
\[
Q(3,4)= ( (1,2,3), (1,2,4), (1,3,4), (2,3,4) ).
\]
Let~$A\in\C^{n\times m }$.   Fix~$k\in\{1,\dots,\min(n,m) \}$.
For~$\alpha  \in Q(k,n),\beta \in Q(k,m)$, let~$A[\alpha|\beta]$ denote the~$ k\times k$ submatrix   obtained by taking the entries of~$A$ in  the rows indexed by~$\alpha$ and the columns indexed by~$\beta$. For example
\[
A[ (2,3)|(1,2)]=\begin{bmatrix}
a_{21} & a_{22}\\
a_{31} &a_{32}
\end{bmatrix} .
\]
The \emph{minor of~$A $ corresponding to~$\alpha ,  \beta$} is
\[
A(\alpha|\beta):=\det (A[\alpha|\beta])  .
\]
For example, if~$m=n$ then~$Q(n,n)$ includes the single element~$\alpha=(1,\dots,n)$, $A[\alpha|\alpha]=A$,
and~$A(\alpha|\alpha)=\det(A)$.
%

Recall that a vector norm $|\cdot |:\mathbb{R}^{n}\to\mathbb{R}_{+}$ induces a matrix norm $||\cdot ||:\mathbb{R}^{n\times n}\to\mathbb{R}_{+}$ defined 
by~$||A||:=\max_{|x|=1} |Ax|$,
and a matrix measure~$\mu (\cdot ):\mathbb{R}^{n\times n}\to\mathbb{R}$   defined by 
    \[
    \mu (A) :=\lim_{\epsilon \downarrow 0} \frac{||I+\epsilon A||-1}{\epsilon}\]
    (see e.g.,~\cite{vid,strom75}). For   the~$L_1$, $L_2$, and~$L_\infty$ norms, there exist closed-form expressions for the induced matrix norms and matrix 
    measures (see Table~\ref{table:matrix_measures}).

\begin{table}[t!]
\centering
{\renewcommand{\arraystretch}{1.5}%
\begin{tabular}{|l|l|l|} 
 \hline
 Vector norm & Induced matrix norm & Induced matrix measure \\ [0.5ex] 
 \hline\hline
 $|x|_1 = \sum_i |x_i|$ & $||A||_1 = \max_j \sum_i |a_{ij}|$ & $\mu_1(A) = \max_j\{a_{jj} + \sum_{i \neq j} |a_{ij}|\}$ \\
 \hline
 $|x|_2 = \sqrt{\sum_i |x_i|^2}$ & $||A||_2 = \sqrt{\lambda_1(A^TA)}$ & $\mu_2(A) = \frac{1}{2}\lambda_1(A + A^T)$ \\ 
 \hline
 $|x|_\infty = \max_i |x_i|$ & $||A||_\infty = \max_i \sum_j |a_{ij}|$ & $\mu_\infty(A) = \max_i\{a_{ii} + \sum_{j \neq i} |a_{ij}|\}$ \\ 
 \hline 
\end{tabular} } 
\caption{  Closed-form expressions for common matrix norms and matrix measures. For a symmetric matrix~$S$, $\lambda_1(S)$ denotes the largest eigenvalue of~$S$. }
\label{table:matrix_measures}
\end{table}
\section{Geometric Motivation}
\begin{figure}
 \begin{center}
  \includegraphics[scale=0.6]{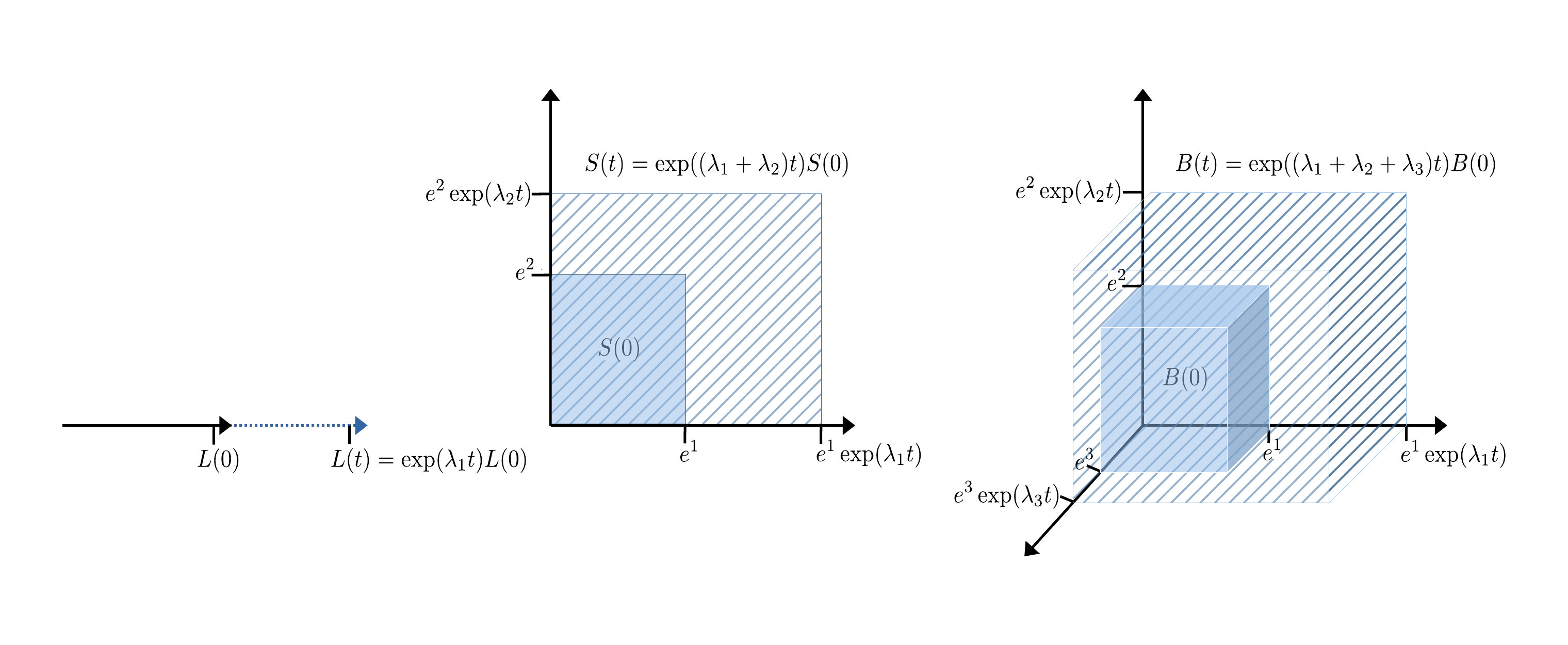}
	\caption{Evolution of lines, areas, and volumes  under the     LTI~\eqref{eq:ltia} with~$\lambda_1>\lambda_2>\lambda_3$.}
	\label{fig:GeometricalEvolution}
\end{center}
\end{figure}

$k$-compound matrices provide information on the evolution of $k$-dimensional parallelotopes   subject  to    a linear time-varying dynamics. To explain this in the  simplest setting,
consider the    LTI system:
\be\label{eq:ltia}
\dot x (t) =\diag(\lambda_1,\lambda_2,\lambda_3) x(t) ,
\ee
with~$\lambda_i\in\R$ and~$x:\R_+ \to\R^3$.
Let~$e^i$, $i=1,2,3$,
denote the~$i$th canonical vector in~$\R^3$.
For~$x(0)=e^i$ we have~$x(t)=\exp(\lambda_i t)x(0)$. Thus, $\exp(\lambda_i t)$ describes the rate of evolution  of the line
between~$0$ and~$e^i$ subject to~\eqref{eq:ltia}.
What about 2D areas? Let~$S_{ij}\subset\R^3 $ denote the square generated by~$e^i$ and~$e^j$, with~$i\not =j$.
Then~$S(t):= x(t,S_{ij} )$ is the rectangle generated by~$\exp(\lambda_i t) e^i$ and~$\exp(\lambda_j t) e^j$, so 
  the area of~$S(t)$ is~$\exp( (\lambda_i + \lambda_j) t)$. Similarly, if~$B_{123} \subset\R^3 $ is the 3D cube generated by~$e^1,e^2$, and~$e^3$ then the volume of~$B(t):=x(t,B_{123})$  
is~$\exp( (\lambda_1+\lambda_2+\lambda_3)t)$ (see Fig.~\ref{fig:GeometricalEvolution}). 
Since~$\exp(At)=\diag ( \exp(\lambda_1 t), \exp(\lambda_2 t), \exp(\lambda_3 t)  )$, this discussion suggests
that it may be useful to have a~$3\times 3$   matrix whose eigenvalues are the sums of any two eigenvalues of~$\exp(At)$, 
and a~$1\times 1$ matrix whose eigenvalue
is the sum of the  three  eigenvalues of~$\exp(At)$.
With this geometric motivation in mind, we   turn to review  the  
  multiplicative and additive compounds of a matrix. For more details and proofs, see e.g.~\cite[Ch.~6]{fiedler_book}\cite{schwarz1970}.

%

\section{Multiplicative Compound}\label{sec:compund_matrices}
%

 The~$k$-multiplicative compound of a matrix~$A$ is a   matrix that collects  all the~$k$-minors of~$A$. 
\begin{Definition}\label{def:multi} 
%
Let $A\in \C^{n\times m}$ 
and fix $k \in \{ 1,\dots,\min  (n,m) \}$. 
The \emph{$k$-multiplicative compound} of~$A$, denoted~$A^{(k)}$, is the~$\binom{n}{k}\times \binom{m}{k}$
matrix that contains all the~$k$-minors of~$A$  ordered lexicographically. 
\end{Definition}

For example, if~$n=m=3$ and~$k=2$ then
\[
A^{(2)}= \begin{bmatrix}
A((12)|(12)) & A((12)|(13)) & A((12)|(23)) \\ 
A((13)|(12)) & A((13)|(13)) & A((13)|(23)) \\ 
A((23)|(12)) & A((23)|(13)) & A((23)|(23)) 
\end{bmatrix}.
\]
In particular, Definition~\ref{def:multi} 
implies that~$A^{(1)}=A$,
and if~$n=m$ then~$A^{(n)}=\det(A)$.
Note also  that by definition~$(A^T)^{(k)}=(A^{(k)})^T$. In particular if~$A$ is symmetric i.e.,~$A=A^T$ 
then~$(A^{(k)})^T =(A^T)^{(k)}= A^{(k)} $, so~$A^{(k)}$ is also symmetric. 
%
 
%
%
%

\subsection{The Cauchy-Binet Formula and its Implications}

The next result, 
known as the  Cauchy-Binet Formula (see e.g.,~\cite{notes_comb_algebra}), justifies the 
 term multiplicative compound.
\begin{Theorem}\label{thm:CB}
Let~$A\in\mathbb{C}^{n\times m}$,  $B\in\mathbb{C}^{m\times p}$. Fix  $k\in \{1,\dots,\min (n,m,p) \}$. Then
\begin{align}\label{eq:AB_MultComp}
    (AB)^{(k)} = A^{(k)} B^{(k)}.
\end{align}
\end{Theorem}
%
%
For~$n=m=p$,  Eq.~\eqref{eq:AB_MultComp}  with~$k=n$ reduces to the familiar formula
\be\label{eq:detsq}
\det(AB)= \det (A) \det (B).
\ee

For the sake of completeness, we include a proof of Thm.~\ref{thm:CB} in the Appendix. 
We now describe some immediate implications of the Cauchy-Binet Formula.

Let~$I_s$ denote the~$s\times s$ identity matrix. 
Definition~\ref{def:multi}  implies that~$I_n^{(k)} = I_r$, where~$r:=\binom{n}{k}$. Hence, if~$A\in\R^{n\times n}$ is non-singular then~$(AA^{-1})^{(k)}=( A^{-1} A)^{(k)}=I_r$ and combining this with~\eqref{eq:AB_MultComp} yields~$(A^{-1})^{(k)} = (A^{(k)})^{-1}$.
In particular, if~$A$ is non-singular then so is~$A^{(k)}$.
Another implication of~\eqref{eq:AB_MultComp} is that if~$T\in\R^{n\times n}$ is non-singular then 
\begin{align}\label{eq:tatm1}
    (TAT^{-1})^{(k)}&= T^{(k)} A^{(k)} (T^{-1} )^{(k)}\nonumber \\
    &=T^{(k)} A^{(k)} (T^{(k)}) ^{-1}.
\end{align}

%
 
%

The  Cauchy-Binet formula also 
yields   a closed-form expression for the spectral properties of
the~$k$-multiplicative compound~$A^{(k)}$ in terms of the  spectral properties of~$A$. 
\begin{Proposition}\label{prop:spect}
For~$A\in\C^{n\times n}$, let~$\lambda_i$,~$i=1,\dots,n$, denote the eigenvalues of~$A$, and let~$v^i$,~$i=1,\dots,n$, denote the eigenvector  corresponding
to~$\lambda_i$. Fix~$k\in\{1,\dots,n\}$.
Then the eigenvalues of~$A^{(k)}$ are all the~$\binom{n}{k}$ products:
\be\label{eq:prodeig}
\bigg\{\prod_{\ell=1}^k \lambda_{i_\ell} \st 1\leq i_1< i_2<\dots< i_k \leq n\bigg\}.
\ee
Furthermore, if
\[
\begin{bmatrix} v^{i_1}&\dots &v^{i_k} \end{bmatrix}^{(k)}.
\]
is not the zero vector then it is 
the eigenvector of~$A^{(k)}$ corresponding to
the eigenvalue~$\prod_{\ell=1}^k \lambda_{i_\ell} $.
\end{Proposition}

To prove Prop.~\ref{prop:spect} note that  
\[
A\begin{bmatrix}
v^{i_1}&\dots& v^{i_k} 
\end{bmatrix} =
\begin{bmatrix}
v^{i_1}&\dots& v^{i_k} \end{bmatrix} \diag(\lambda_{i_1} ,\dots,\lambda_{i_k}),
\]
and 
applying the Cauchy-Binet formula gives 
\begin{align*}
    A^{(k)} \begin{bmatrix}
v^{i_1}&\dots& v^{i_k} 
\end{bmatrix}^{(k)} & =
\begin{bmatrix} v^{i_1}&\dots &v^{i_k} \end{bmatrix}^{(k)}
(\diag(\lambda_{i_1} ,\dots,\lambda_{i_k}))^{(k)}\\
&= \left( \prod_{\ell=1}^k \lambda_{i_\ell} \right) \begin{bmatrix} v^{i_1}&\dots &v^{i_k} \end{bmatrix}^{(k)}.
\end{align*}

\begin{Example}
Suppose that~$n=3$ 
and~$A$ is upper-triangular: 
$A=\begin{bmatrix}
a_{11} & a_{12} & a_{13}\\
0 & a_{22} & a_{23}\\
 0 &  0  & a_{ 33}
\end{bmatrix}  .
$
Then the eigenvalues of~$A$ are~$\lambda_i=a_{ii}$, $i=1,2,3$.
Assume that~$c:=   a_{22}-a_{11} \not = 0 $. Then the  
  first two corresponding eigenvectors 
  are~$v^1=\begin{bmatrix}1&0&0
\end{bmatrix}^T$, 
$v^2=\begin{bmatrix}  
a_{12} &
c&0
\end{bmatrix}^T$,
and~$v^1$ and~$v^2$ are linearly  independent. Now  
\begin{align*}
\begin{bmatrix}
v^1 & v^2
\end{bmatrix}^{(2)}=
\begin{bmatrix}
 1&  a_{12} \\
 0&  c  \\
 0&0
\end{bmatrix}^{(2)}= c \begin{bmatrix} 1 \\0\\ 0 
\end{bmatrix}.
\end{align*}
A direct  calculation gives
\[
A^{(2)}= \begin{bmatrix}
a_{11}   a_{22} &  a_{11} a_{23}&  a_{12} a_{23}-a_{13} a_{22} \\
0& a_{11} a_{33}& a_{12} a_{33}\\
0& 0& a_{22} a_{33}
 \end{bmatrix} ,
\]
so the eigenvalues of~$A^{(2)}$
are of the  form~\eqref{eq:prodeig}, and the first eigenvector of~$A^{(2)}$ is indeed~$\begin{bmatrix}
v^1 & v^2
\end{bmatrix}^{(2)}$. 
\end{Example}

The next result describes another  simple and useful application of the multiplicative compound.


\begin{Proposition}\label{prop:gendet}
Suppose that~$X,Y\in\C^{n\times k}$, with~$k\leq n$. Then
\be\label{eq:gen_det}
\det(Y^T X)=(Y ^{(k)})^T X ^{(k)}.
\ee
\end{Proposition}
Note that~$X^{(k)}$ (and~$Y^{(k)}$) has dimensions~$\binom{n}{k}\times\binom{k}{k}$, i.e. it is a column vector, so the right-hand side of~\eqref{eq:gen_det} is the inner product of two vectors. In the particular case where~$k=n$, Eq.~\eqref{eq:gen_det} reduces to 
$\det(Y^T X)= \det(Y)\det(X)$. However,~\eqref{eq:gen_det} holds also for \emph{non-square} matrices.

\begin{IEEEproof}
Using the fact that~$Y^TX\in\R^{k\times k}$ and the  Cauchy-Binet formula yields
 \begin{align*} 
     \det(Y^TX) 
     &=(Y^T X) ^{(k)}\nonumber\\
     &=(Y^T) ^{(k)} X ^{(k)}\nonumber\\
     &=(Y ^{(k)})^T X ^{(k)},
 \end{align*}
and this completes the proof.
\end{IEEEproof}

Prop.~\ref{prop:gendet} has found many applications. For example, Ref.~\cite{city14026} applied it to determine the sensitivity of the natural modes of an electrical circuit  to modifications in the circuit  elements and topology. As we will see in the next section,
Prop.~\ref{prop:gendet} also implies that the 
 $k$-multiplicative compound can be used
 to  describe the  volume of~$k$-dimensional parallelotopes.

 Another implication of the Cauchy-Binet Formula is that
 certain  sign properties of the minors of a matrix are preserved under matrix multiplication. To explain this, we 
 require the following definition.
 \begin{Definition}\label{def:tnr}\cite{total_book}
 Let~$A\in \R^{n\times m}$, and
  fix~$r\in\{1,\dots,\min(n,m) \}$. Then~$A$ 
 is called \emph{totally non-negative of order~$r$} ($TN_r$) if every minor of size~$\leq r$ of~$A$ is non-negative.
 $A$ 
 is called \emph{totally positive of order~$r$}~($TP_r$) if every minor of size~$\leq r$ of~$A$ is positive. 
 \end{Definition}
 
 Such matrices admit    sign variations diminishing properties that have important applications to dynamical systems (see Section~\ref{sec:posi} below).  In general, $TN_r$ and~$TP_r$ are not preserved under natural matrix operations.
 \begin{Example}\label{exa:not_sum} 
 Consider 
    the matrices~$A=\begin{bmatrix} 1&1 &1 \\ 1&1 &1 \\1&1 &1\end{bmatrix}$ and~$B=I_3$. A direct calculation of all the minors shows that both these matrices are~$TN_3$.
    However, 
 $A+B=\begin{bmatrix}
 2&1&1\\ 1&2&1\\ 1& 1&2
 \end{bmatrix}$
 admits  a minor of order two that is negative, as~$\det(\begin{bmatrix}
 1&1 \\2 &1
 \end{bmatrix})=-1$, so~$A+B$   is not~$TN_3$ (and not even~$TN_2$). 
 \end{Example}
 However, the  Cauchy-Binet Formula immediately  implies that the \emph{product}
 of two~$TN_r$ [$TP_r$] matrices is~$TN_r$ [$TP_r$].

%

\subsection{Multiplicative Compounds   and the Volume of Parallelotopes}\label{sec:volumes}
We review an important interpretation of the $k$-multiplicative compound     based on the presentation in~\cite[Chapter~IX]{Gantmacher_vol1} and~\cite{GOVER201028}.
For a vector~$v\in\R^{n}$, let~$|v|_2 := \sqrt{v^{T} v}$   denote the $L_{2}$ norm of $v$.
Fix~$k\in\{1,\dots,n\}$ and vectors~$x^1,\dots,x^k \in \R^n$. The parallelotope generated by these vectors (and the zero vertex) is the set:
\[
P:=\bigg\{ \sum_{i=1}^k r_i x^i\st r_i \in [0,1] \bigg\}.
\]
Note that this implies that~$0\in\R^n$ is a vertex of~$P$. 
This  can always be assured  by a simple translation.
We can also interpret~$P$
 as the image of the unit~$k$-cube
 under the matrix
 \[
 X:=\begin{bmatrix}
 x^1&\dots& x^k
 \end{bmatrix}\in\R^{n\times k} .
 \]
The~\emph{Gram matrix} associated with~$x^1,\dots,x^k$
is the $k\times k$ symmetric matrix:
\begin{align}\label{eq:defgram}
    G ( x^1,\dots,x^k): &=X^T X \nonumber\\
    &=
    \begin{bmatrix}
    (x^1) ^T x^1 & (x^1) ^Tx^2 & \hdots & (x^1) ^T x^k \\
     &   \vdots \\
     (x^k) ^T x^1 & (x^k) ^Tx^2 & \hdots & (x^k) ^T x^k 
    \end{bmatrix}.
\end{align}
For example, for~$k=1$ we have~$G(x^1)=|x^1|_2^2$, and for~$k=2$ we have
\[ 
    G ( x^1, x^2) =
    \begin{bmatrix}
    |x^1|_2^2 & (x^1) ^Tx^2   \\
     (x^2) ^T x^1 & |x^2|_2^2   
    \end{bmatrix}.
\]

It follows from~\eqref{eq:defgram}
  that for any~$s\in\R^k$ we have
\[
|\sum_{i=1}^k s_i x^i  |_2^2 = s^T G(x^1,\dots,x^k) s,
\]
so~$G(x^1,\dots,x^k)$ is non-negative definite, 
and   it is positive-definite iff~$x^1,\dots,x^k$ are linearly independent.  

\begin{figure}
 \begin{center}
  \includegraphics[scale=0.2]{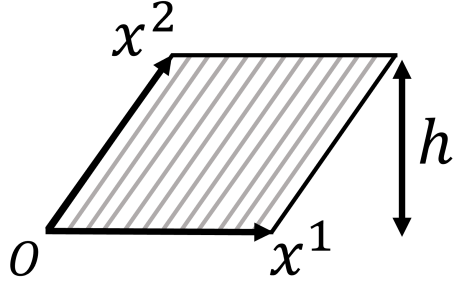}
	\caption{The volume of the parallelotope $P(x^{1},x^{2})$ is $\text{vol} ( P(x^{1}) ) h = |x^{1}|_2  h$.}
	\label{fig:Parallelogram_k=2}
\end{center}
\end{figure}

\begin{figure}
 \begin{center}
  \includegraphics[scale=0.2]{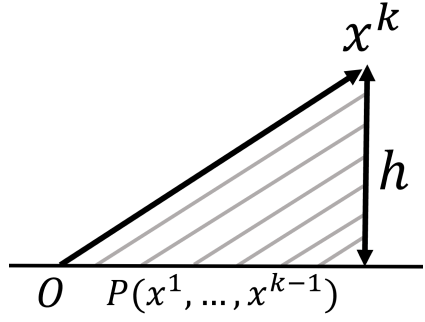}
	\caption{The volume of the parallelotope $P(x^{1},\hdots,x^{k})$ is $\text{vol} ( P ( x^{1},\hdots,x^{k-1} ) ) h$.}
	\label{fig:Parallelogram_k}
\end{center}
\end{figure}

The 
  volume of~$P(x^1,\dots,x^k)$ is
denoted~$\vol(P(x^1,\dots,x^k))$. For~$k=1$, the parallelotope is just the line connecting the origin  and~$x^1$, so~$\vol(P(x^1))=|x^1|_2$. 
For~$k=2$, the parallelotope~$P(x^1,x^2)$ is depicted in Fig~\ref{fig:Parallelogram_k=2}. In general,  $\vol(P(x^1,\dots,x^k))$  is defined in  a recursive manner.
\begin{Definition}\label{def:vol_parall}
 For~$k=1$,~$\vol(P(x^1)):=|x^1|_2$. For any~$k>1$,
\be\label{eq:def_vol}
\vol(P(x^1,\dots,x^k)):=\vol(P(x^1,\dots, x^{k-1}) ) h,
\ee
where~$P(x^1,\dots, x^{k-1})$
 is the~$(k-1)$-dimensional  ``base'' of~$P(x^1,\dots,x^ k)$,
and the altitude~$h$ is the distance from~$x^k$ to the base (see Fig.~\ref{fig:Parallelogram_k}).
\end{Definition}


The next result relates  $\text{Vol} (P)$ to the $k$-multiplicative 
compound of the matrix~$X$. 
\begin{Proposition}\label{prop:vol_via_gram}
The volume of~$P(x^1,\dots,x^k)$ satisfies: 
\be\label{eq:vol_via_compound}
\vol(P(x^1,\dots,x^k))=|X ^{(k)}|_2.
 \ee
\end{Proposition}
Note that since~$X\in\R^{n\times k }$, $X ^{(k)}$ is a column vector. 
 
\begin{IEEEproof}
We first prove a simple
algebraic expression for the volume in terms of the Gram matrix, namely, 
\be\label{eq:volpeqn}
\vol(P(x^1,\dots,x^k))=\sqrt{\det(G(x^1,\dots,x^k))}. 
\ee
To prove this, let~$y$ denote the ``foot'' of the altitude~$h$, that is, 
\be\label{eq:y=sumx}
y=\sum_{i=1}^{k-1} r_ix^i, \text{ with } r_i \in \R, 
\ee
and
\be\label{eq:linsys1}
       0=( x^k-y)^T x^j, \quad j=1,\dots,k-1 .
\ee
Then the altitude~$h$ satisfies
\begin{align}\label{eq:linsys2}
h^2&=|x^k-y|_2^2 \nonumber\\
   &=(x^k-y)^T(x^k-y)\nonumber\\
   &=(x^k-y)^T x^k.
\end{align}
Eqs.~\eqref{eq:linsys1} and~\eqref{eq:linsys2}  yield the linear system
\be\label{eq:lin_sys}
\begin{bmatrix} 
(x^1)^T x^1 & (x^2)^T x^1 & \hdots &       (x^{k-1})^T x^1 &0 \\
&&\vdots \\
(x^1)^T x^{k-1} & (x^2)^T x^{k-1} & \hdots &       (x^{k-1})^T x^{k-1}&0 \\
(x^1)^T x^{k } & (x^2)^T x^{k } & \hdots &       (x^{k-1})^T x^{k }&1 
\end{bmatrix}
\begin{bmatrix} r_1\\\vdots\\r_{k-1} \\h ^2\end{bmatrix} = \begin{bmatrix} (x^k)^T x^1 \\\vdots\\ (x^k)^T x^{k-1} \\ (x^k)^T x^k  \end{bmatrix}.
\ee

If~$x^1,\dots,x^k$ are linearly  dependent, 
that is,~$x^k$ lies in the~$(k-1)$-dimensional subspace 
that contains  the ``base'' of the parallelotope then~$h=0$, so~\eqref{eq:def_vol} gives~$\vol(P(x^1,\dots,x^k))=0$ 
and~\eqref{eq:volpeqn} gives the same value. 

If~$x^1,\dots,x^k$ are linearly  independent then 
 applying Cramer's rule~\cite[Ch.~0]{matrx_ana} to~\eqref{eq:lin_sys} 
 gives
\begin{align*}
h^2&=\frac{\det\begin{pmatrix} 
(x^1)^T x^1 & (x^2)^T x^1 & \hdots &       (x^{k-1})^T x^1 &  (x^{k })^T x^1 \\
&&\vdots \\
(x^1)^T x^{k-1} & (x^2)^T x^{k-1} & \hdots &       (x^{k-1})^T x^{k-1}&  (x^{k })^T x^{k-1} \\
(x^1)^T x^{k } & (x^2)^T x^{k } & \hdots &       (x^{k-1})^T x^{k }&  (x^{k })^T x^ k
\end{pmatrix} }
{\det \begin{pmatrix} 
(x^1)^T x^1 & (x^2)^T x^1 & \hdots &       (x^{k-1})^T x^1 &0 \\
&&\vdots \\
(x^1)^T x^{k-1} & (x^2)^T x^{k-1} & \hdots &       (x^{k-1})^T x^{k-1}&0 \\
(x^1)^T x^{k } & (x^2)^T x^{k } & \hdots &       (x^{k-1})^T x^{k }&1 
\end{pmatrix} },
\end{align*}
that is,
\begin{align*}
h^2  = \frac{\det(G(x^1,\dots,x^{k}))}{\det(G(x^1,\dots,x^{k-1}))},
\end{align*}
and combining this with~\eqref{eq:def_vol}  and Def.~\ref{def:vol_parall}  proves~\eqref{eq:volpeqn}.
Now using the fact that
    $\det(G(x^1,\dots,x^k))=
     \det(X^TX)$,
     and Prop.~\ref{prop:gendet} 
     proves~\eqref{eq:vol_via_compound}.
\end{IEEEproof}

Note that in the special case where~$k=n$, the matrix~$X$ is a square matrix, and~\eqref{eq:vol_via_compound}
gives
$(\vol(P(x^1,\dots,x^k)) )    ^2 = 
 (\det(X))^2$,
 Thus, for~$k=n$ we recover
 the well-known formula
    \[
    \vol(P(x^1,\dots,x^n))=|\det(\begin{bmatrix}
    x^1&\dots & x^n
    \end{bmatrix})|.
    \]

\subsection{Multiplicative Compounds and the Sign Variations Diminishing Property}\label{subsec:SVDP}
One reason for the usefulness of compound matrices in systems and control theory is that they allow
to track the evolution of volumes along the solutions of differential equations. 
Another important application of compound matrices is in the context of sign variation diminishing properties. 
We review this topic in the classical setting of totally positive matrices~\cite{gk_book,pinkus,total_book}.
Unfortunately, this field suffers from non-uniform notation. We follow the notation in~\cite{total_book}.

 \subsubsection{Totally Positive Matrices}
\begin{Definition}
A matrix~$A\in\R^{n\times m}$ is called totally positive~(TP)
if all its minors are   positive.
\end{Definition}
Equivalently, $A\in\R^{n\times m}$ is  TP  if it is~$TP_r$ with~$r=\min(n,m)$. 
Another equivalent definition is that for any~$k\in\{1,\dots, \min( n,m)\}$ 
the $k$-multiplicative compound matrix~$A^{(k)}$ has positive entries. For example, for~$A=\begin{bmatrix}
2&1&1\\1&3&4
\end{bmatrix}$, we have~$A^{(1)}=A\gg 0$ and~$A^{(2)}=\begin{bmatrix}
5&7&1
\end{bmatrix} \gg 0$, so~$A$ is~TP. 
If~$n=m$, i.e.,~$A$ is square, then~$A$ is TP if and only 
if~$A^{(k)}$ maps~$\R^{\binom{n}{k}}_+\setminus\{0\}$ to~$\R^{\binom{n}{k}}_{++}$ for every~$k\in\{1,\dots,n\}.$

  TP matrices satisfy  a beautiful sign variation diminishing property. 
For a vector~$x\in\R^n\setminus\{0\}$, let~$s^-(x) $ denote the number of sign variations in~$x$ after deleting all  its
zero entries. For example,~$s^-(\begin{bmatrix}-1&0&0&2&-3\end{bmatrix}^T)=2$. We define~$s^-(0):=0$. For a vector~$x\in\R^n $,  let~$s^+(x)$
denote the maximal possible number of sign variations in~$x$ after setting  every zero entry in~$x$ to either~$-1$ or~$+1$. For
example,~$s^+(\begin{bmatrix}-1&0&0&2&-3\end{bmatrix}^T)=4$. These definitions imply that
\be\label{eq:sminussplus}
0\leq s^-(x)\leq s^+(x)\leq n-1,   \text{ for all } x\in\R^n.
\ee

There is  a useful duality relation between~$s^-$ and~$s ^+$ that is straightforward to prove. 
Let~$D_{\pm}:=\diag(1,-1,1,\dots,(-1)^{n-1})$. Then
\[
s^-(x)+s^+(D_{\pm} x)=n-1,\text{ for any } x\in\R^n. 
\]

The next result is useful when studying the number of sign variations in a vector that is the limit of  a sequence of vectors.
\begin{Proposition}\cite[Chapter 3]{pinkus}\label{prop:limitschanges}
Suppose that~$x^1,x^2,\dots  $
is a sequence of vectors in~$ \R^n$ that converges to a limit~$x^*:=\lim_{i \to \infty} x^i$. 
Then
\begin{align}\label{eq:schlimit}
    s^-(x^*)&\leq \liminf_{i\to \infty } s^-(x^i)  \leq  \limsup_{i\to \infty } s^+(x^i)
  \leq  s^+(x^*).  
\end{align}
\end{Proposition}
Intuitively speaking, the number of sign changes in the limit vector~$x^*$ can change only
if there exists an index~$k$ such that~$x^*_k=0$. Eq.~\eqref{eq:schlimit}
follows from the fact that  zero entries  are ignored in~$s^-$, but may lead to an increase in~$s^+$. 
For example,  for~$x^i=\begin{bmatrix}-1& 2^{-i}&-1\end{bmatrix}^T$  we have~$s^-(x^i)= 2$ for any~$i$, $s^-(x^*) =0$, and~$s^+(x^*)=2$.

The next result describes the sign variation diminishing property~(SVDP) of~TP matrices.
\begin{Theorem}\label{thm:svdptp}\cite{Pinkus1996}
Suppose that~$A\in\R^{n\times n}$ is~TP.
Then  for any~$x\in\R^n\setminus\{0\}$, we have 
\be\label{eq:SVDP}
s^+(Ax)\leq s^-(x).
\ee
Furthermore,
if
\be\label{eq:eqsigns}
s^+(Ax) = s^- (x)
\ee
then the sign of the first [last]
component of~$Ax$ (if zero, the sign given in determining $s^+(Ax)$)
agrees with the sign of the first  [last] nonzero component of~$x$.
\end{Theorem}

In other words, multiplying a vector by a TP matrix can never increase the number of sign variations between~$s^-(x)$ and~$s^+(Ax)$,
and if the number of sign changes remains the same then in some sense~$x$ and~$Ax$   have the   same  ``orientation''.

\begin{Example}
Consider~$A=\begin{bmatrix}
1&b\\c&d
\end{bmatrix}$,
with
\be\label{eq:bcd}
b,c,d>0 \text{ and }  d>bc.
\ee
Then~$A$ is TP.
Fix~$x\in\R^2\setminus\{0\}$. Then one of the following cases holds.  

\noindent \emph{Case 1.}  Suppose that~$s^-(x)=s^+(x)=0$. We may assume with out loss of generality~(wlog) that~$x_1,x_2>0$.
Then
\[
y:=Ax=\begin{bmatrix}
x_1+b x_2\\cx_1+dx_2
\end{bmatrix}
\]
satisfies~$y_1,y_2>0$, so~$s^+(y)=0 =s^-(x)$, 
 and  the sign of the first [last]
component of~$y$ 
agrees with the sign of the first  [last] nonzero component of~$x$.

\noindent \emph{Case 2.}  Suppose that~$s^-(x)=0$ and~$s^+(x)=1$. We may assume wlog that~$x_1=0,x_2=1$.
Then
\[
y =\begin{bmatrix}
 b  \\d
\end{bmatrix}
\]
satisfies~$y_1,y_2>0$, so~$s^+(y)=0 =s^-(x)$, 
 and   the sign of the first [last]
component of~$y$ 
agrees with the sign of the first  [last] nonzero component of~$x$.

\noindent \emph{Case 3.}  Suppose that~$s^-(x)=s^+(x)=1$. We may assume wlog that~$x_1=1,x_2<0$.
Then
\[
y =\begin{bmatrix}
1+b x_2\\c+dx_2
\end{bmatrix}
\]
and~\eqref{eq:bcd} implies that~$y_2< c y_1$. If~$y_1<0$ then~$y_2<0$, 
so~$s^+(y)<  s^-(x)$.  
If~$y_1=0$ then~$y_2<0$,
so~$s^+(y)=s^-(x)$, and the sign     given to~$y_1$ in determining~$s^+(y)$ is plus. If~$y_1>0 $ then either: (1)~$y_2>0$ and then~$s^+(y)<s^-(x)$;
or (2)~$y_2= 0$
and then~$s^+(y)=s^-(x)$ and the sign     given to~$y_2$ in determining~$s^+(y)$ is minus; or~(3)~$y_2< 0$
and then~$s^+(y)=s^-(x)$.

Thus, we see that in each  case the assertions in Thm.~\ref{thm:svdptp} holds. 
\end{Example}

As we will see in Section~\ref{sec:posi}
below, 
the SVDP has important implications in the asymptotic analysis 
of linear and non-linear  dynamical systems~\cite{margaliot2019revisiting}.

 \subsubsection{Recognition of Totally Positive Matrices}
Since the minors of a matrix are not independent, verifying that a matrix is~TP does not require checking that all minors are positive. This fact will play an important role in Section~\ref{sec:hankel}, so we   review one result  on the recognition of~TP matrices. For more details, see~\cite{total_book}.

\begin{Definition}
   An index set~$\alpha\in Q(k,n)$  is called \emph{contiguous} if it consists of only consecutive numbers i.e., $\alpha_1=p, \alpha_2=p+1,\dots,\alpha_k=p+k-1$ for some integer~$p$. If~$\alpha,\beta \in Q(k,n)$ are two contiguous index sets   then the corresponding submatrix~$A[\alpha\mid\beta]$ is called a~\emph{contiguous submatrix of}~$A$, and the~$k$-minor~$A(\alpha\mid\beta)$ is called a~\emph{contiguous minor}. A contiguous minor is called~\emph{initial} if at least one of the two sets~$\alpha$ or~$\beta$ is~$\{1,2,\hdots,k\}$. 
\end{Definition}

The next result shows that   if  certain initial and contiguous minors are positive  then   all minors are positive. 
\begin{Proposition}\label{prop:Fallat}~\cite{Fallat_2017}
Let~$A\in\R^{m\times n}$. If all the initial minors of~$A$ 
  up to   order~$k-1$ are positive, and all its contiguous~$k$-minors are positive, then
     $A $ is~$TP_k$.
\end{Proposition}
  
\begin{Example}\label{exp:falla}
    Consider the matrix
$
        A=\begin{bmatrix}
        3&1&2\\2&1&3\\1&3&10
        \end{bmatrix}.
$
We apply  Prop.~\ref{prop:Fallat} to verify that~$A^{(2)}\gg0$, i.e., that    all the~$2$-minors of~$A$ are positive. The initial minors up to  order 1 are the   the entries in the first row and first column of~$A$, and these  are positive. The  contiguous~$2$-minors
are~$A((1,2)\mid(1,2)), A((1,2)\mid(2,3)), A((2,3)\mid(1,2))$ and~$A((2,3)\mid(2,3))$, which are~$1,1,5$, and~$1$, respectively, so we 
conclude that~$A^{(2)}\gg0$.
\end{Example}

\section{Additive Compound}\label{sec:add_comp}
The matrix~$\exp(At)$ describes the evolution of lines 
under the LTI dynamics~$\dot x=Ax$. We will see below that, more generally,  the evolution of~$k$-dimensional parallelotopes is described by~$(\exp(At))^{(k)}$. This naturally leads to the following  question:
what is the derivative with respect to~(w.r.t.) time of~$(\exp(At))^{(k)}$? To address this, we require the $k$-additive compound of~$A$.

\begin{Definition}\label{def:add_comp}
%
Let~$A\in\mathbb{C}^{n\times n}$.
The \emph{$k$-additive compound}
matrix of~$A$ is the~$\binom{n}{k}\times \binom{n}{k}$ matrix defined by:
\begin{align}\label{eq:A^[k]:=ddeps}
    A^{[k]} := \frac{d}{d\epsilon} (I_n+\epsilon A)^{(k)} |_{\epsilon=0} . 
\end{align}
%
\end{Definition}
The derivative here is well-defined, as every entry of~$(I_n+\epsilon A)^{(k)}$ is a polynomial in~$\epsilon$.
Note that this definition implies that
\be\label{eq:akasexp}
A^{[k]} =  \frac{d}{d\epsilon} (\exp(A\epsilon ))^{(k)} |_{\epsilon=0} ,
\ee
and also that
\begin{align}\label{eq:(I+epsilonA)^k}
    (I_n+\epsilon A)^{(k)} = I_r + \epsilon A^{[k]} + o(\epsilon),
\end{align}
%
where~$r:=\binom{n}{k}$. 
In other words,~$A^{[k]}$ is the coefficient of the first-order term in the Taylor expansion of~$(I_n+\epsilon A)^{(k)}$.  The definition also  implies that
\begin{align}\label{eq:dexp/dt=A^[k]exp}
    \frac{d}{dt}(\exp(At))^{(k)}=A^{[k]}(\exp(At))^{(k)}
\end{align}
(see  the more general result in  Prop.~\ref{prop:odeforphik} below). 
 Thus,~$(\exp(At))^{(k)}$ satisfies an LTI with  the matrix~$A^{[k]}$. Note that for~$k=1$, Eq.~\eqref{eq:dexp/dt=A^[k]exp} reduces to~$ \frac{d}{dt}\exp(At)=A\exp(At)$, whereas for~$k=n$ it becomes
 $\frac{d}{dt}\det(\exp(At))=\tr(A)\det(\exp(At))$.

Eq.~\eqref{eq:A^[k]:=ddeps} and the Cauchy-Binet formula can be used
to determine   how~$A^{[k]}$ changes under a similarity transformation of~$A$.
If~$T$ is non-singular then
 \begin{align}\label{eq:sim_trans}
 (TAT^{-1})^{[k]}&=  \frac{d}{d\epsilon} (I_n+\epsilon TAT^{-1})^{(k)} |_{\epsilon=0}\nonumber \\
 &= \frac{d}{d\epsilon} (T(I_n+\epsilon A)T^{-1})^{(k)} |_{\epsilon=0}\nonumber \\
 &=T^{(k)} A^{[k]} (T^{(k)})^{-1}.
 \end{align}

The next result describes the spectral properties of the additive compound. Its proof follows from combining 
Prop.~\ref{prop:spect} and~\eqref{eq:(I+epsilonA)^k}.
\begin{Proposition}\label{prop:spect_add}
For~$A\in\C^{n\times n}$, let~$\lambda_i$,~$i=1,\dots,n$, denote the eigenvalues of~$A$, and let~$v^i$,~$i=1,\dots,n$, denote the eigenvector  corresponding
to~$\lambda_i$. Fix~$k\in\{1,\dots,n\}$.
Then the eigenvalues of~$A^{[k]}$ are all the~$\binom{n}{k}$ sums:
\[
\bigg\{\sum_{\ell=1}^k \lambda_{i_\ell} \st 1\leq i_1< i_2<\dots< i_k \leq n\bigg\}.
\]
Furthermore, if
\[
\begin{bmatrix} v^{i_1}&\dots &v^{i_k} \end{bmatrix}^{(k)}
\]
is not the zero vector then it is 
the eigenvector of~$A^{[k]}$ corresponding the eigenvalue~$\sum_{\ell=1}^k \lambda_{i_\ell} $.
\end{Proposition}
In particular, if~$A$ is positive-definite (so it is symmetric and all its  eigenvalues are real  and positive) then~$A^{[k]}$ is symmetric and all its eigenvalues are real and positive, so~$A^{[k]}$ is positive-definite. 

Another important  implication  of the definitions above is that for any~$A,B\in\C^{n\times n}$ we have
\be
\label{eq:addi}
(A+B)^{[k]}=A^{[k]}+B^{[k]}.
\ee
This justifies the term additive compound. 
Moreover, the mapping~$A\to A^{[k]}$ is linear.
%
To prove~\eqref{eq:addi}, note that~\eqref{eq:(I+epsilonA)^k}
gives
\begin{align*}
I_r+\epsilon (A+B)^{[k]}&=     
      (I_n+\epsilon (A+B))^{(k)}+o(\epsilon)\\
   &=\left (  (I_n+\epsilon A)    (I_n+\epsilon B) \right) ^{(k)}+o(\epsilon)\\
&=  (I_n+\epsilon A)^{(k)}    (I_n+\epsilon B)   ^{(k)}+o(\epsilon)\\&=
\left (  I_r+\epsilon A^{[k]} \right)\left (  I_r+\epsilon B^{[k]} \right)+o(\epsilon)\\
&=I_r+\epsilon (A^{[k]}+B^{[k]} )+o(\epsilon)
\end{align*}
and  using the continuity of the mapping~$A\to A^{[k]}$  implies~\eqref{eq:addi}.  
%
%

The next  result gives a useful explicit formula for~$A^{[k]}$ in terms of the entries~$a_{ij}$ of~$A$.
  Recall that any entry of~$A^{(k)}$ is a minor~$A(\alpha|\beta)$.
Thus, it is natural to index the entries of~$A^{(k)}$  and~$A^{[k]}$  using~$\alpha,\beta \in Q(k,n)$. 
%
\begin{Proposition}\label{prop:Explicit_A_k}
%
Fix~$\alpha,\beta \in Q(k,n)$ and let~$\alpha=\{i_1,\dots,i_k\}$ and~$\beta=\{j_1,\dots,j_k\}$. Then the entry of~$A^{[k]}$ corresponding to~$(\alpha,\beta)$ is
equal to:
\begin{enumerate}
    \item $\sum_{\ell=1}^{k} a_{ i_{\ell} i_{\ell} }$, if  $i_{\ell} = j_{\ell}$  for all $\ell \in \{ 1,\hdots,k \}$;
    \item $(-1)^{\ell +m}
    a_{i_{\ell} j_{m}} $, 
      if all the indices  in $ \alpha  $  and $ \beta$   agree,  except  for   
   a single index $ i_{\ell} \ne j_m$; and
    \item $0$, otherwise.
\end{enumerate}
%
\end{Proposition}
 Note that the   first case in the proposition corresponds to the diagonal entries of~$A^{[k]}$. Also,  the proposition
implies in particular that~$A^{[n]}= \sum_{\ell=1}^{n} a_{\ell \ell}= \tr(A)$.

We prove  Prop.~\ref{prop:Explicit_A_k}  for the case~$k=2$. 
The proof when~$k>2$ is similar.  
Let $B:= (I_n+\epsilon A)^{(2)}$. Denote
\begin{align*}
    \delta_{pq} := 
    \begin{cases}
    1, & \text{if } p=q, \\
    0, & \text{otherwise}.
    \end{cases}
\end{align*}
Then for any $1\le i_1 < i_2 \le n$ and $1 \le j_1 < j_2 \le n$,
\begin{align*}
    B( (i_1,i_2) | (j_1,j_2) ) &= 
    b_{i_1 j_1} b_{i_2 j_2} - b_{i_1 j_2} b_{i_2 j_1} \\
    &= ( \delta_{i_1 j_1} + \epsilon a_{i_1 j_1} )
       ( \delta_{i_2 j_2} + \epsilon a_{i_2 j_2} ) 
        - 
       ( \delta_{i_1 j_2} + \epsilon a_{i_1 j_2} )
       ( \delta_{i_2 j_1} + \epsilon a_{i_2 j_1} ) \\
       &= \epsilon ( \delta_{i_1 j_1} a_{i_2 j_2}
       +\delta_{i_2 j_2} a_{i_1 j_1}
       -\delta_{i_1 j_2} a_{i_2 j_1}
       -\delta_{i_2 j_1} a_{i_1 j_2})
       +c + o(\epsilon),
\end{align*}
where $c$  does not depend on $\epsilon$.
Applying~\eqref{eq:(I+epsilonA)^k} yields 
\begin{align*}
    A^{[2]} ( \alpha|\beta ) =
    \delta_{i_1 j_1} a_{i_2 j_2} +
    \delta_{i_2 j_2} a_{i_1 j_1} -
    \delta_{i_1 j_2} a_{i_2 j_1} -
    \delta_{i_2 j_1} a_{i_1 j_2},
\end{align*}
and this agrees with the expressions in Prop.~\ref{prop:Explicit_A_k}.
%
\begin{Example}\label{ex:A^[3]}
%
For~$A\in\R^{4\times 4}$ and~$k=3$,
Prop.~\ref{prop:Explicit_A_k} yields
%
%
\[
A^{[3]}  =
\begin{blockarray}{ccccc}
(1,2,3) & (1,2,4) & (1,3,4) & (2,3,4) \\
\begin{block}{[cccc]c}
  a_{11}+a_{22}+a_{33} & a_{34} & -a_{24} & a_{14} & (1,2,3) \\
  a_{43} & a_{11}+a_{22}+a_{44} & a_{23} & -a_{13} & (1,2,4) \\
  -a_{42} & a_{32} & a_{11}+a_{33}+a_{44} & a_{12} & (1,3,4) \\
  a_{41} & -a_{31} & a_{21} & a_{22}+a_{33}+a_{44} & (2,3,4) \\
\end{block}
\end{blockarray},
 \]
where the indexes~$\alpha \in Q(3,4)$ [$ \beta \in Q(3,4)$]
are marked on  right-hand side  
[above]  of  the matrix.
For example, the entry in the second row and fourth column
of $A^{[3]}$ corresponds to 
$(\alpha | \beta) = ((1,2,4 ) |  ( 2,3,4 ) )$.
As~$\alpha$ and~$\beta$ agree in all indices except
for the  index
$i_{ 1} =1$ and $j_{ 2} =3$,
this  entry is equal to~$  (-1)^{1+2} a_{13} = -a_{13}$.
%
\end{Example}
\begin{Example}\label{ex:A23}
%
For~$A\in\R^{3\times 3}$ and~$k=2$,
Prop.~\ref{prop:Explicit_A_k} yields
%
%
\be\label{eq:add23}
A^{[2]}  =\begin{blockarray}{cccc}
(1,2) & (1,3) & (2,3)  \\
\begin{block}{[ccc]c}
  a_{11}+a_{22}  & a_{23} & -a_{13} & (1,2) \\
  a_{32} & a_{11}+a_{33} & a_{12} & (1,3) \\
  -a_{31} & a_{21} & a_{22}+a_{33} & (2,3) \\
\end{block}
\end{blockarray}.
\ee
For example, the 
entry in the second row and third  column
of $A^{[3]}$ corresponds to 
$(\alpha | \beta) = ((1,3 ) |  ( 2,3 ) )$.
As~$\alpha$ and~$\beta$ agree in all indices except
for the  index
$i_{ 1} =1$ and $j_{ 1} =2$,
this  entry is equal to~$  (-1)^{1+1} a_{12} = a_{12}$.
%
\end{Example}

  The next section describes applications of compound matrices for dynamical systems described by~ODEs. 
  For more details and proofs, see~\cite{schwarz1970,muldowney1990}.

%
\section{Compound Matrices and ODEs}\label{sec:odes}
%
Fix a time interval~$[\tau_0,\tau_1] $. Let~$A:[\tau_0,\tau_1] \to \R^{n\times n} $ be a continuous matrix function, and consider the
linear time-varying~(LTV) system:
\be\label{Eq:ltv}
\dot x(t)=A(t)x(t),\quad x(\tau_0)=x_0. 
\ee
The solution at time~$t$ is~$x(t)=\Phi(t ,\tau_0)x_0$,
where~$\Phi(t,\tau_0 )$ is the solution at time~$t$ of the matrix differential 
equation
\be\label{eq:phidot}
\frac{d}{ds} \Phi(s)=A(s) \Phi(s), 
\quad \Phi(\tau_0)=I_n. 
\ee
Fix~$k\in\{1,\dots,n\}$ and let~$r:=\binom{n}{k}$.
A natural question is: how do the~$k$-order minors of~$\Phi(t)$ evolve in time?
\subsection{$k$-Compound System}
The next result provides an elegant  formula for the evolution of~$\Phi^{(k)}(t):= (\Phi(t))^{(k)}$.
\begin{Proposition}\label{prop:odeforphik}
If~$\Phi$ satisfies~\eqref{eq:phidot} then
\be\label{eq:expat}
\frac{d}{ds} \Phi^{(k)}(s)=A^{[k]}(s) \Phi^{(k)}(s),\quad \Phi^{(k)}(\tau_0)=I_r,
\ee
where~$A^{[k]}(s):= (A(s)) ^ {[k]} $.
\end{Proposition}
\begin{IEEEproof}
For~$\varepsilon>0$, we  have
\begin{align*}
\Phi^{(k)}(s+\varepsilon)& = \left ( \Phi(s)+\varepsilon A(s) \Phi(s) +o(\varepsilon) \right ) ^{(k)}\\
&=\left ( (I_n+\varepsilon A(s) )  \Phi(s)    \right ) ^{(k)} +o(\varepsilon) \\
&= (I_n+\varepsilon A(s) )   ^{(k)}  \Phi^{(k)}(s) +o(\varepsilon) \\
&=( I_r + \varepsilon A^{[k]} (s) )\Phi^{(k)}(s) +o(\varepsilon) ,
\end{align*}
so 
\begin{align*}
\frac{ \Phi^{(k)}(s+\varepsilon) -\Phi^{(k)}(s ) }{\varepsilon} & =  
  A^{[k]} (s)  \Phi^{(k)}(s) +\frac{o(\varepsilon)}{\varepsilon} ,
\end{align*}
and this  completes the proof. 
\end{IEEEproof}

Eq.~\eqref{eq:expat} is the \emph{$k$-compound system} of~\eqref{eq:phidot}, and is the basis of all the 
$k$-generalizations of dynamical systems
described in this paper.

The dynamics of the  $k$-compound system 
implies that~$k$-minors of~$\Phi$ also satisfy an~LTV.
In particular, if~$A(t)\equiv A$  and~$\tau_0=0$
then~$\Phi(t)=\exp( At) $ so~$\Phi^{(k)}(t)=(\exp(At))^{(k)}$, and~\eqref{eq:expat} gives
\be\label{eq:expatti}
(\exp(At))^{(k)} = 
\exp(A ^{[k]} t).
\ee
This identity has interesting implications. 
For example, recall that~$\exp(Bt)\geq 0$ (where the inequality is component-wise) for all~$t\geq 0$ iff~$B$ is Metzler. Eq.~\eqref{eq:expatti} implies that~$(\exp(At))^{(k)}\geq 0$ for all~$t\geq 0$ iff~$A^{[k]}$ is Metzler. This is the basis for the notion of~$k$-positive systems described in Section~\ref{sec:posi} below.

Note also that for the special case~$k=n$, Prop.~\ref{prop:odeforphik}  yields 
 \[
 \frac{d}{dt} \det(\Phi(t))
 = \tr(A(t)) \det(\Phi(t)),
 \]
which is the
Abel-Jacobi-Liouville identity (see, e.g.,~\cite[p. 152]{chicone_2006}).

Prop.~\ref{prop:odeforphik} implies that under the LTV dynamics~\eqref{Eq:ltv}, $k$-dimensional parallelotopes evolve  according to the dynamics~\eqref{eq:expat}. 
Indeed, 
consider $k$ initial conditions~$a^{1},\hdots,a^{k} \in\R^{n}$. Under the LTV dynamics, and assuming for simplicity that~$t_0=0$,  the 
  corresponding solutions satisfy~$x  (t,a^i) =\Phi (t) a^{ i} $.
Let~$X(t):=\begin{bmatrix}
x(t,a^1)& \hdots& x(t,a^k)
\end{bmatrix}$. Then the 
  volume   of the parallelotope~$P(x (t,a^1),\dots,x  (t,a^k)) $  is~$|X^{(k)}(t) |_2$.
Now,
\begin{align}\label{eq:xktdyn}
    X^{(k)}(t)
    &=
    \begin{bmatrix}
    \Phi (t) a^{1} & \hdots & \Phi (t) a^{k}
    \end{bmatrix}^{(k)}\nonumber \\
    &=
    \left (
    \Phi (t)
    \begin{bmatrix}
    a^{1} & \hdots & a^{k}
    \end{bmatrix} 
    \right ) ^{(k)}\nonumber\\
    &=
    \Phi^{(k)} (t)
        \begin{bmatrix}
    a^{1} & \hdots & a^{k}
    \end{bmatrix}^{(k)}\nonumber\\
    &=
    \Phi^{(k)} (t)
         X^{(k)}( 0), 
\end{align}
where the third equation follows from the Cauchy Binet formula. Using~\eqref{eq:expat} gives
\begin{align}\label{eq:aaddeq}
    \frac{d}{dt} \left ( X^{(k)} (t) \right )
    &= 
    \frac{d}{dt} \left ( \Phi^{(k)} (t) X^{(k)} (0) \right )\nonumber\\
    &= A^{[k]} (t) \Phi^{(k)} (t) X^{(k)} (0) \nonumber\\
    &= A^{[k]} (t) X^{(k)} (t).
\end{align}
%
%
Note that this    implies in particular that if~$X^{(k)}(0)\not =0 $ then~$X^{(k)}(t)\not =0$ for all~$t $. In other words, if~$a^1,\dots,a^k$ are linearly independent then~$x(t,a^1),\dots,x(t,a^k)$ are linearly independent for all~$t$.

The compounds can also be used in the analysis  of non-linear dynamical systems. Consider the  time-varying system
\be\label{eq:nonlin}
\dot x (t) = f(t,x).
\ee
For the sake of simplicity, we  assume that the initial time is zero, and that the system admits a convex and compact state-space~$\Omega$.
We also assume that~$f\in C^1$. The Jacobian of the vector field~$f$ is~$J(t,x):=\frac{\partial}{\partial x}f(t,x)$.
Compound matrices can be used to analyze~\eqref{eq:nonlin} by using an LTV called the  variational equation associated  with~\eqref{eq:nonlin}.
To define it, 
fix~$a,b \in \Omega$.
Let~$z(t):=x(t,a)-x(t,b)$, and for~$s\in[0,1]$, 
 let~$\gamma(s):=s x(t,a)+(1-s)x(t,b)$, i.e. the line connecting~$x(t,a)$ and~$x(t,b)$.
 Then
\begin{align*} 
   \dot z(t)&=f(t,x(t,a))-f(t,x(t,b))  \\
	       &=\int_0^1  \frac{\partial }{\partial s}   f(t,\gamma(s))     \dif s  ,
\end{align*}
and this gives the variational equation:
\be\label{eq:var_eqn}
\dot z(t)=A^{ab}(t)z(t),
\ee
where
\begin{equation}\label{eq:at_int}
    A^{ab}(t):=\int_0^1  J(t,\gamma(s))     \dif s.  
\end{equation}
Note that the variational equation~\eqref{eq:var_eqn} is an LTV. 

In the next section, we use the results above to  describe a general principle   for deriving useful
``$k$-generalizations'' of  important 
classes of non-linear 
dynamical systems including cooperative systems~\cite{hlsmith},  contracting systems~\cite{sontag_cotraction_tutorial,LOHMILLER1998683}, 
and diagonally stable systems~\cite{kaszkurewicz2012matrix}. 
\section{$k$-Generalizations of Dynamical Systems}\label{sec:k_generalizations}
%
Suppose that  the LTV~\eqref{Eq:ltv}
  satisfies a specific \emph{property} guaranteed  by a related condition on~$A(t)$. For example,  \emph{property} may be that the LTV is positive (i.e.~$A(t)$ is Metzler for all~$t$)  or that the~LTV is contracting (guaranteed if~$\mu(A(t))\leq-\eta<0$ for all~$t$, where~$\mu:\R^{n\times n}\to\R$ is some matrix measure).
  Fix~$k\in\{1,\dots,n\}$.  We say that
the~LTV satisfies~\emph{$k$-property} if the~$k$-compound system 
satisfies the related property. For example, the LTV is~\emph{$k$-positive} if~$A^{[k]}(t)$ is Metzler for all~$t$;  the LTV is~\emph{$k$-contracting} if~$\mu(A^{[k]}(t))\leq-\eta<0$ for all~$t$, and so on. 

This generalization principle  was first suggested in~\cite{comp_barshalom}. It makes sense for two reasons. First, when~$k=1$, the~$k$-compound system  reduces to the original system, so 
\emph{$k$-property} is clearly a generalization of
\emph{property}. Second, the $k$-compound system  has a clear geometric meaning:  it describes the evolution of~$k$-dimensional parallelotopes along the dynamics. 
 Also,  the $k$-compound system   describes the evolution of $k$ minors,  and 
 the sign of these minors is important in establishing an~SVDP.

The same principle  can be applied to  the nonlinear system~\eqref{eq:nonlin}
using the variational equation~\eqref{eq:var_eqn}.
For example, if~$\mu(J(t,x))\leq -\eta<0$ for all~$t\geq 0$ and all~$x\in\Omega$ then~\eqref{eq:nonlin} is contracting: the distance between any two solutions (in the norm that induced~$\mu$) decays at an exponential rate. If we replace this by the condition~$\mu(J^{[k]}(t,x))\leq -\eta<0$ for all~$t\geq 0$ and all~$x\in\Omega$
(i.e., the same condition but now for the $k$-compound system)  then~\eqref{eq:nonlin} is called~$k$-contracting. Roughly speaking, 
this means that the  volume of $k$-dimensional 
parallelotopes  decays to zero exponentially along the flow of the nonlinear system.
We now turn to describe several  such~$k$-generalizations in more detail.  
%
\section{$k$-Contracting Systems}\label{sec:kcont}
%
The term $k$-contracting systems was  introduced in~\cite{kordercont}. For~$k=1$ these reduce to contacting systems. 
This   generalization of contracting systems   is motivated in part by the seminal work of Muldowney~\cite{muldowney1990}
who analyzed  nonlinear systems that, in the new terminology, are~$2$-contracting (see also the unpublished preprint~\cite{weak_manchester} 
for some preliminary ideas). Muldowney derived several interesting results for   $2$-contracting systems. For example, every bounded trajectory of a \emph{time-invariant}  $2$-contracting system
converges to an equilibrium (but, unlike in the case of contracting systems, the equilibrium is not necessarily unique). These results have found many applications in 
models from epidemiology, see e.g.,~\cite{SEIR_LI_MULD1995}. Such models typically have two equilibrium points corresponding to the 
disease free   and the endemic steady-states. 
The existence of more than a single equilibrium point implies that
the system is not contracting (that is, not $1$-contracting) w.r.t. any norm. 

In this tutorial, we focus on explaining  $k$-contraction in LTVs. The case of nonlinear systems follows by applying the similar  ideas to the variational  equation~\eqref{eq:var_eqn}, see~\cite{kordercont}. 

\begin{Definition}\label{def:conlin} \cite{kordercont}
Fix~$k\in \{1,\dots,n\} $. 
The LTV~\eqref{Eq:ltv}  (with initial time~$\tau_0=0$)
is  called \emph{$k$-contracting} if there exist~$\eta > 0$ and a vector norm $|\cdot|$ such that for
any~$a^1, \dotsm, a^k \in \R^n$,   the mapping $X(\cdot):\R_+\to  \R^{n \times k}$ defined by~$
X(t) : = \begin{bmatrix} x(t, a^1) & \dots & x(t, a^k)
\end{bmatrix}
$
satisfies
\be\label{eq:contdef}
| X^{(k)}(t)|
 \leq \exp(- \eta t)
 | X^{(k)}(0) |, \text{ for all }t \geq 0.
\ee
\end{Definition}

In other words,   the volume of  any  $k$-parallelotope  exponentially   decays to zero   under the dynamics. Another interpretation 
is that the initial condition
\[
X^{(k)} (0)   = \begin{bmatrix}  a^1 & \dots &  a^k
\end{bmatrix}^{(k)}
\]
is ``forgotten''.

For~$k=1$, Eq.~\eqref{eq:contdef}    reduces to the requirement that~$|x(t,a)|\leq \exp(- \eta t)|a|$ for any~$a\in\R^n$, and since the dynamics is linear this implies 
that~$|x(t,a)-x(t,b)|\leq \exp(- \eta t)|a-b|$ for any~$a,b\in\R^n$
i.e., standard  contraction.   

Eq.~\eqref{eq:xktdyn}
implies that~$k$-contraction is equivalent to 
\[
||\Phi^{(k)}(t) || \leq \exp(-\eta t) || \Phi^{(k)}(0)|| = \exp(-\eta t) ,\text{ for all } t\geq 0. 
\]

Using Coppel's inequality~\cite{coppel1965stability} 
 and~\eqref{eq:aaddeq} provides  a  simple sufficient condition for   $k$-contraction in terms of the $k$ additive compound of~$A$.  
\begin{Proposition}
If there exist~$\eta>0$ and a matrix measure~$\mu:\R^{n\times n}\to\R$ such that 
\be\label{eq:kconinfi} 
\mu(A^{[k]}(t))\leq-\eta<0, \text{ for all } t \geq 0,
\ee
then~\eqref{Eq:ltv} is~$k$-contracting.
\end{Proposition}

\begin{figure}
 \begin{center}
  \includegraphics[scale=0.8]{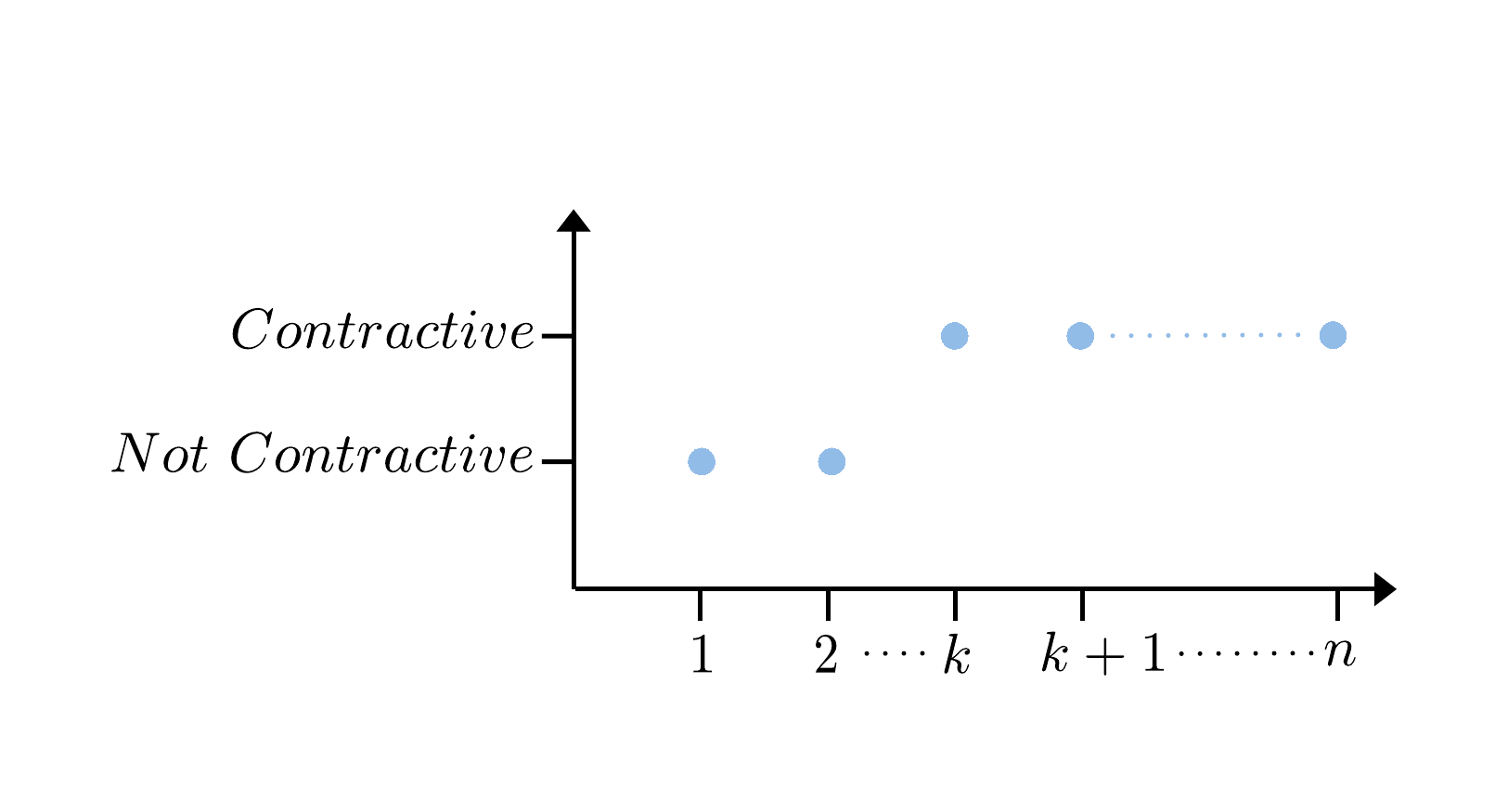}
	\caption{$k$-contraction implies~$\ell$-contraction for any integer~$\ell\geq k$.}
	\label{fig:k_cont}
\end{center}
\end{figure}

Note that for~$k=1$ condition~\eqref{eq:kconinfi}  reduces to the standard infinitesimal condition for contraction~\cite{sontag_cotraction_tutorial}.
For~$k=n$,  
condition~\eqref{eq:kconinfi} becomes~$\tr(A(t))\leq -\eta<0 $ for all $t \geq 0$.
It was shown in~\cite{kordercont} that if~\eqref{eq:kconinfi}   holds for some~$L_p$ norm, with~$p\in\{1,2,\infty\}$  then
  for any integer~$\ell \geq k$ we have
$
	 \mu ( J^{[\ell]}(t ) )\leq -\eta<0$
for all~$t\geq 0$, so  the   system is also  $\ell$-contracting 
w.r.t. the~$L_p$ norm (see Fig.~\ref{fig:k_cont}).

\begin{Example}
Consider the LTI system~$\dot{x}=Ax$ with~$A=\begin{bmatrix}0&c\\-c&0
\end{bmatrix}$,~$c\in\R$. Since~$\frac{d}{dt}(x_1^2(t)+x_2^2(t))\equiv 0$, the solution for any initial condition~$x(0)$ is a circle with radius~$r:=\sqrt{x_1^2(t)+x_2^2(t)} \equiv\sqrt{x_1^2(0)+x_2^2(0)}$.
Fix~$a^1,a^2\in\R^2$ and consider~$X(t):=\begin{bmatrix}x(t,a^1)&x(t,a^2)\end{bmatrix}$.   Eq.~\eqref{eq:aaddeq}   gives
\begin{align*}
   \frac{d}{dt}\left(X^{(2)}(t)\right)&=A^{[2]}X^{(2)}(t)\\&=\trace(A)X^{(2)}(t)\\&=0.
\end{align*}
This shows that the
LTI is``on the verge'' of being~$2$-contracting. This makes sense because the matrix~$A$ represents a purely rotational dynamics, so  the area of  the  parallelotope generated by~$x(t,a^1)$ and~$x(t,a^2)$ is a constant of time for any value~$c$.
\end{Example}
 
For the~$L_p$ norms, with~$p\in\{1,2,\infty\}$, 
  condition~\eqref{eq:kconinfi}     is easy to check, as 
  combining Prop.~\ref{prop:Explicit_A_k} with the expressions in Table~\ref{table:matrix_measures} gives~\cite{muldowney1990}:
\begin{align*} 
\mu_1(A^{[k]}) &= 
 \max_{\alpha \in Q(k,n) }  \Big( \sum_{p=1}^k a_{\alpha_p,\alpha_p}  
 + \sum_{\substack{j \notin \alpha }}(|a_{j,\alpha_1}| + \cdots + |a_{j,\alpha_k}|) \Big) ,\nonumber \\
\mu_2(A^{[k]}) &= \sum_{i=1}^k \lambda_i (  {A + A^T}   )/2 , \\
\mu_{\infty}(A^{[k]}) &=
\max_{\alpha \in Q(k,n) } \Big(  \sum_{p=1}^k a_{\alpha_p,\alpha_p} + \sum_{\substack{j \notin\alpha }}(|a_{\alpha_1,j}| + \cdots + |a_{\alpha_k,j}|) \Big) ,
\end{align*}
where~$\lambda_i(A+A^T)$, $i=1,\dots,n$, denote the eigenvalues of the symmetric matrix~$A+A^T$ 
  ordered such that~$\lambda_1\geq \dots\geq \lambda_n$.

\begin{figure}
 \begin{center}
  \includegraphics[scale=0.5]{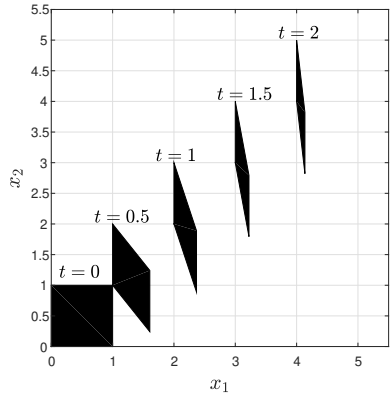}
	\caption{Evolution of  the unit square in Example~\ref{exa:squares}.}
	\label{fig:squares}
\end{center}
\end{figure}

\begin{Example}\label{exa:squares}
	Consider the LTV~\eqref{Eq:ltv} with~$n=2$, $\tau_0=0$,  and
$
	A(t) = \begin{bmatrix}
	 -1 & 0 \\
	-2\cos( t) &0  
	\end{bmatrix}.
$
The corresponding  transition matrix is:
$
						  	\Phi(t)= \begin{bmatrix} \exp(-t)&0 \\ 
										                -1+\exp(-t)(\cos(t)-\sin(t))  &1\end{bmatrix}.
	$	This implies that the LTV is uniformly stable, and that for any~$x(0)\in\R^2$ we have
	\be\label{eq:limxt2}
	\lim _{t\to \infty} x(t,x(0))=\begin{bmatrix}  0  \\ x_2(0)-x_1(0)  \end{bmatrix} .
	\ee
	The LTV 
 is not contracting w.r.t.   any norm, as   it admits   
 more than a single equilibrium. However,~$ A^{[2]}(t) =\tr(A(t)) \equiv -1$, so the 
	system is~$2$-contracting.  Let~$S_0\subset\R^2$ denote the unit square, and let~$S(t):=x(t,S_0 )$, that is,
	the evolution
	at time~$t$
	of the unit square   under the dynamics. 
	Fig.~\ref{fig:squares} 
	depicts~$S(t) +2t$  for several values of~$t$,
	where the shift by~$2t$ is for  the sake of clarity.
	It may be seen that the area of~$S(t)$ decays with~$t$, and that~$S(t)$ converges to a line.
\end{Example}
Note that convergence to sets that are more general than an equilibrium point is a  common requirement  in applications. For example, in consensus algorithms~\cite{eger_consens} the requirement is that the state~$x(t)$ converges to the subspace~$\spanop(\begin{bmatrix}
1&\dots& 1
\end{bmatrix}^T)$.

The next result describes another 
interesting implication of~$k$-contraction of an LTV.
\begin{Proposition}\label{prop:eqi2} \cite{muldowney1990}
\label{prop:sub}  
	Suppose that the LTV~\eqref{Eq:ltv} is uniformly stable. Fix~$k\in\{1,\dots,n\}$.
The following   conditions are equivalent: 
	\begin{enumerate}[(a)]
	\item \label{cond:nk1} The LTV  
	admits an~$(n-k+1)$-dimensional linear subspace 
	$\mathcal{X} \subseteq \R^n$ such that
	\begin{equation} \label{eq:a0x}
	\lim_{t \to \infty }x(t, x_0) = 0 \text{ for any } x_0 \in \mathcal{X} .
	\end{equation}
	\item \label{cond:sysmkop} Every solution of  
	\begin{equation} \label{eq:kinsys}
	\dot{y}(t) = A^{[k]}(t)y(t) , 
	\end{equation}
	with~$y(\cdot)\in \R^{\binom{n}{k}}$,
	converges to the origin.
\end{enumerate}
\end{Proposition}

For~$k=1$ the equivalence of conditions~$(a)$ and~$(b)$ is obvious. For~$k=n$, condition~$(a)$ 
becomes the existence of a one-dimensional subspace~$\mathcal{X} \subseteq \R^n$ such that any solution emanating from~$\mathcal{X} 
$ converges to the origin, whereas~$(b)$
requires that every solution of the scalar system~$\dot y(t)=\tr(A(t))y(t)$ converges to the origin.

Note that if~$\mu(A^{[k]}(t))\leq-\eta<0$ for all~$t\geq 0$ then condition~$(b)$  holds, and thus 
condition~$(a)$ holds.

\begin{Example}
	Consider the special case where the  LTV~\eqref{Eq:ltv} reduces to an LTI~$\dot x(t)=Ax(t)$. The uniform stability requirement implies that the real part of every eigenvalue of~$A$ is non-positive.  
	Condition~$(b)$ is equivalent to the requirement that the sum of every~$k$ eigenvalues of~$A$ has  a negative real part. 
	This implies that at least~$(n-k+1) $ eigenvalues of~$A$ have a negative real part, and thus 
	condition~$(a)$ holds.
\end{Example}

\begin{Example}
	Consider again the LTV in Example~\ref{exa:squares}. Here~$n=2$ and the LTV is uniformly stable and~$k$-contracting for~$k=2$,  so Prop.~\ref{prop:eqi2}
	implies that
	the LTV  
	admits a one-dimensional linear subspace 
	$\mathcal{X} \subseteq \R^2$ such that~\eqref{eq:a0x} holds. Indeed, it follows from~\eqref{eq:limxt2}
	that~$\spanop(\begin{bmatrix}
	1&1
	\end{bmatrix}^T)$ is this subspace. 
\end{Example}
 
As noted above,
the definition of~$k$-contractivity for non-linear systems follows by applying the same approach to the variational equation which is an~LTV. We refer to~\cite{kordercont} for the details. Muldowney and his colleagues~\cite{muldowney1990,li1995}   proved that~\emph{time-invariant}  $2$-contracting systems have a ``well-ordered''
asymptotic behaviour, and this has been used to derive a global asymptotic
analysis of important non-linear dynamical   
models from epidemiology (see, e.g.~\cite{SEIR_LI_MULD1995}). 
A recent paper~\cite{searial12contracting} extended some of these results to systems that are not necessarily~$2$-contracting, but
can be represented as the serial interconnections of~$k$-contracting systems, with~$k\in\{1,2\}$. 
Ref.~\cite{ofir2021sufficient} studied the series connection of two systems, and derived a sufficient condition for $k$-contraction of the overall system. This is based on a new formula for the $k$-compounds of a block-diagonal matrix. 

 The notion of~$k$-contracting systems is based on the relation between $k$-compounds and the volume of $k$-dimensional parallelotopes. However, the $k$-minors of a matrix have another important application
 related to  the sign variation diminishing  property.
 Using this allows  to introduce~$k$-generalizations of  the important classes of linear positive systems and non-linear cooperative systems.

\section{$k$-Positive Systems}\label{sec:posi}
 Positive linear systems~\cite{farina2000} and cooperative non-linear systems~\cite{hlsmith}  
 are characterized  by the fact that every  state-variable takes   non-negative values. This is the case in many real-world systems where the state-variables represent quantities such as densities in queues, probabilities, concentration of molecules, etc. An important property  of such  systems is that many  system and control problems scale well with the system size~\cite{rantzer_valcher,tanaka2011}.

Ref.~\cite{Eyal_k_posi} introduced the notions of~$k$-positive and~$k$-cooperative systems. The LTV~\eqref{Eq:ltv} is called \emph{$k$-positive} if~$A^{[k]}(t)$ is
Metzler   for all~$t\in[\tau_0,\tau_1]$. In other words, we require the standard condition for positivity, 
but on the $k$-compound system. 
For~$k=1$ this reduces to requiring
that~$A(t)$ is Metzler for all~$t\in[\tau_0,\tau_1]$. In this case the system is positive i.e. the flow maps~$\R^n_+$ to~$\R^n_+$ (and also~$\R^n_-:=-\R^n_+$ to~$\R^n_-$)~\cite{farina2000}.
In other words, the flow maps the set of vectors with zero sign variations to itself. 
The flow of~$k$-positive systems maps the set of vectors with up to~$(k-1)$ sign variations to itself. To explain this, we use  the definitions and results introduced in Section~\ref{subsec:SVDP}.

For any $k\in \{1,\dots,n\}$,
define the sets~\cite{fusco1988}\cite[p.~71]{krasno_posi_opera}:
\begin{align}\label{eq:Pksets}
P_{-}^{k} &:= \{ z\in\mathbb{R}^{n}\st \; s^{-} (z) \le k-1 \},\nonumber \\
P_{+}^{k} &:= \{ z\in \mathbb{R}^{n}\st \; s^{+} (z) \le k-1 \}.
\end{align}
In other words,
these are the sets of all vectors with up to~$(k-1)$ sign variations. 
Then~$P_{-}^{k}$ is closed, and it can be shown that~$P_{+}^{k}=\Int(P_{-}^{k})$. 
For example,
\begin{align*}
    P_{-}^{1} = \mathbb{R}_{+}^{n} \cup \mathbb{R}_{-}^{n}, \;\;\; 
    P_{+}^{1} = \Int(\mathbb{R}_{+}^{n}) \cup \Int( \mathbb{R}_{-}^{n}). 
\end{align*}

Since multiplying a vector by a positive (or a negative) constant does not change the number of sign variations in the vector, the sets~$P^k_-$ and~$P^k_+$ are cones. However, they are not convex cones. For example, for~$n=2$, $x=\begin{bmatrix}
2&1
\end{bmatrix}^T$ and~$y=\begin{bmatrix}
-1&-2
\end{bmatrix}^T $, we have~$s^-(x)=s^-(y)=s^+(x)=s^+(y)=0$, but~$s^-(\frac{x+y}{2})= s^+(\frac{x+y}{2})=1$.
For an analysis  of  the geometric structure of these sets, see~\cite{Eyal_k_posi}.

\begin{Definition}
The LTV~\eqref{Eq:ltv} is called \emph{$k$-positive} on the  interval~$[\tau_0,\tau_1]$ if for any~$\tau_0<t_0  < \tau_1$,
\[
x_0\in P^k_- \implies x(t,t_0,x_0) \in P^k_- \text{ for all }   t_0 \leq t < \tau_1 ,
\]
and is called \emph{strongly $k$-positive} if
\[
x_0\in P^{k}_{-}\setminus \{0\} \implies x(t,t_0,x_0) \in P^k_+ \text{ for all }  \tau_0 < t < \tau_1. 
\]
\end{Definition}
In other words, the sets of up to~$(k-1)$ sign variations are invariant sets of the dynamics.  

Recall the important SVDP of TP matrices~\eqref{eq:SVDP}, namely,    multiplying a vector by a TP matrix can never increase the number of sign variations. For our purposes, we need a more specialized result.  A matrix~$A\in\R^{n\times m}$ is called \emph{sign-regular of order~$k$} (denoted~$SR_k$) if its minors of order~$k$ are all non-positive or all non-negative. It is called  \emph{strictly 
sign-regular of order~$k$} (denoted~$SSR_k$)
if   its minors of order~$k$ are    all positive or all   negative. 
In this case, we use~$\epsilon_k$ to denote the signature of the~$k$-minors, that is,~$\epsilon_k=1$ [$\epsilon_k=-1$] if all the~$k$-minors are positive [negative].
%
For example,
$
    A = \begin{bmatrix}
    1 & 2 \\
    3 & 4
    \end{bmatrix}
$
 is~$SSR_1$ with~$\epsilon_1=1$ as every entry
 (i.e., every $1$-minor)   is positive, and it is~$SSR_2$
 with~$\epsilon_2=-1$
 since~$\det(A)<0$, and this is the only $2$-minor
 of~$A$.
%
%
\begin{Proposition}\label{thm:BenAvraham_SSRk}
%
(see e.g.~\cite{CTPDS})
Let $A\in\mathbb{R}^{n\times n}$ be a non-singular matrix. Fix $k\in \{1,\dots,n\}$. Then the following two conditions are equivalent:
\begin{enumerate}
    \item For any $x\in\mathbb{R}^{n}  $ with
    $s^{-}(x) \le k-1$, we have~$s^{-} (Ax) \le k-1$.
    \item $A$ is   $SR_k$.
\end{enumerate}
Also,  the following two conditions are equivalent:
\begin{enumerate}[a)]
    \item For any $x\in\mathbb{R}^{n} \setminus \{ 0 \}$ with
    $s^{-}(x) \le k-1$, we have~$s^{+} (Ax) \le k-1$.
    \item $A$ is  $SSR_k$.
\end{enumerate}
\end{Proposition}

\begin{Example}
Consider the non-singular matrix
$
    A=
    \begin{bmatrix}
    3&2&-1 \\3&5&-1\\3&5&0
    \end{bmatrix} .
$
Then
$
    A^{(2)}=
    \begin{bmatrix}
   9&0&3  \\ 9&3&5 \\ 0&3&5
    \end{bmatrix}
$, so~$A$ is $SR_2$. 
Fix\be\label{eq:repo} 
x\in\R^3 \setminus\{0\} \text{ with }s^-(x)\leq 1.
\ee
Prop.~\ref{thm:BenAvraham_SSRk} implies that~$s^-(Ax)\leq 1$. We  verify this directly. Seeking a contradiction, assume that~$s^- (Ax)>1$ i.e.,~$s^-(Ax)=2$. 
Then we may assume, wlog, that~$y:=Ax$ satisfies~$y_1,y_3 > 0$ and~$y_2<  0$, that is,
\begin{align*}
    3x_1+2x_2-x_3 &> 0,\\
   3 x_1+5 x_2-x_3 & < 0,\\
    3 x_1+ 5 x_2   & >  0.
\end{align*}
The first two equations give~$x_2<  0$, and the 
  last two equations give~$x_3 >  0$.
  Now the first  equation implies that~$x_1 > 0$,  and 
  this contradicts~\eqref{eq:repo}.
\end{Example}

Note that the statements in Prop.~\ref{thm:BenAvraham_SSRk} do  not directly compare~$s^-(A x )$ and~$ s^-(x)$.
However, the results in this proposition  imply the following. 
%
\begin{Corollary}\label{coro:ssr_for_all_j}
Let $A\in\mathbb{R}^{n\times n}$ be a non-singular matrix. Fix $k\in \{1,\dots,n\}$. Then the following two conditions are equivalent:
\begin{enumerate}
    \item \label{item:first_cor} For any $x\in\mathbb{R}^{n}  $ with
    $s^{-}(x) \le k-1$, we have~$s^{-} (Ax) \leq s^{-}(x)$.
    \item \label{item:sec_cor}  $A$ is    $SR_j$ for every~$j\in\{1,\dots,k\}$.
\end{enumerate}
Also,  the following two conditions are equivalent:
\begin{enumerate}[a)]
    \item For any $x\in\mathbb{R}^{n} \setminus \{ 0 \}$ with
    $s^{-}(x) \le k-1$, we have~$s^{+} (Ax) \le s^{-}(x) $.
    \item $A$ is  $SSR_j$ for every~$j\in\{1,\dots,k\}$.
\end{enumerate}
\end{Corollary}

\begin{IEEEproof}
Suppose that~$A$ is  $SR_j$ for every~$j\in\{1,\dots,k\}$. Fix~$x\in\R^n$ with~$s^-(x)\leq k-1$. 
Then there exists a~$j \in\{1,\dots,k\}$ such that~$s^-(x) = j-1$, and Prop.~\ref{thm:BenAvraham_SSRk} implies that~$s^-(Ax)\leq j-1$, so~$  s^-(Ax) \leq s^-(x) $.

Now assume that there exists some~$j\in\{1,\dots,k\}$ such that~$A$ is \emph{not}  sign-regular of order~$j$.
Then Proposition~\ref{thm:BenAvraham_SSRk} implies that there exists~$x\in\R^n$ such that~$s^-(x) \leq j-1$ and~$s^-(Ax)>j-1$, so~$s^-(Ax)>s^-(x) $. 
This proves the equivalence of the first two assertions. 
The remainder of the proof is similar, and thus omitted. 
\end{IEEEproof}

\begin{Remark}\label{rem:thnr}

Note that Corollary~\ref{coro:ssr_for_all_j} implies the following. 
Let~$A\in\R^{n\times n}$ be non-singular. 
If~$A $ is~$TN_k$    
then  $s^{-}(x) \le k-1$ implies that~$s^{-} (Ax) \leq s^{-}(x)$.
If~$A $ is~$TP_k$    
then  for any~$x\not =0$ with~$s^{-}(x) \le k-1$ we have~$s^{+} (Ax) \leq s^{-}(x)$.
\end{Remark}

Using these tools allows to characterize~$k$-positive LTVs.
\begin{Theorem}\label{thm:k_pois_ltv} \cite{Eyal_k_posi}
The LTV~\eqref{Eq:ltv}
is~$k$-positive on~$[\tau_0,\tau_1]$ iff~$A^{[k]}(t)$ is Metzler  
for all~$t\in (a,b)$.
It is strongly $k$-positive  on~$[\tau_0,\tau_1] $
iff~$A^{[k]}(t)$ is Metzler for all~$t\in(\tau_0,\tau_1 )$, 
and~$A^{[k]}(t)$ is irreducible  for all~$t \in(\tau_0,\tau_1 )$  except, perhaps, 
at isolated time points.
\end{Theorem}

The proof is simple. Consider for example  the second assertion in the theorem. The  Metzler and irreducibility assumptions 
imply that the matrix differential system~\eqref{eq:expat} is a positive linear system, and furthermore,   all the entries of~$\Phi^{(k)}(t,t_0)$ are positive for all~$t>t_0$ (recall that the initial condition is~$\Phi(t_0)=I\geq 0$). In other words,~$\Phi(t,t_0)$ is~$SSR_k$ for all~$t>t_0$. Since~$x(t,t_0, x(t_0))=\Phi(t,t_0)x(t_0)$,  applying Prop.~\ref{thm:BenAvraham_SSRk} completes the proof. 

This line of reasoning  demonstrates the  general   principle in Section~\ref{sec:k_generalizations}, namely, given conditions on~$A^{[k]}$ we can   apply standard tools from dynamical  systems theory  to the $k$-compound system~\eqref{eq:expat}, and deduce results on the behaviour of the solution~$x(t)$ of the original system~\eqref{Eq:ltv}. 

\subsection{Sign Conditions for $k$-positivity}
A natural question is: when is~$A^{[k]}$ a Metzler matrix? This can be answered using
Prop.~\ref{prop:Explicit_A_k} in terms of   sign pattern conditions on the entries~$a_{ij}$ of~$A$. This is useful, as in   fields like
chemistry and systems biology,   exact values of various parameters are typically unknown, but their signs may be inferred from various properties of the system~\cite{sontag_near_2007}.
\begin{Proposition}\label{prop:sign_pattern_for_k_posi} \cite{Eyal_k_posi}
Let~$A \in\R^{n\times n}$ with~$n\geq 3$.
Then
\begin{enumerate}
\item \label{case:nminus1}
$A^{[ n-1]}$ is Metzler iff
  $a_{ij}\geq 0$ for all~$i,j$ with~$i-j$ odd, and~$a_{ij}\leq 0$ for all~$i\not = j$ with~$i-j$ even;
\item \label{case:kodd}
for any  odd~$k  $ in the range~$1<k<n-1$,  
$A^{[k]}$ is Metzler   
 iff
 $a_{1n},a_{n1}\geq 0$, $a_{ij}\geq 0$ for all~$|i-j|=1$, and~$a_{ij}=0$ for all~$1<|i-j|<n-1$;
\item \label{case:keven} 
for
any  even~$k  $ in the range~$1<k<n-1$,  $A^{[k]}$ is Metzler   
 iff
 $a_{1n},a_{n1}\leq 0$, $a_{ij}\geq 0$ for all~$|i-j|=1$, and~$a_{ij}=0$ for all~$1<|i-j|<n-1$.
\end{enumerate}
\end{Proposition}

\begin{Example}
For~$n=5$, the corresponding sign patterns are as follows. The matrix~$A^{[4]}$ is Metzler iff
\be\label{eq:a45}
A=\begin{bmatrix}
*& + & -& + &- \\
+ & *& + & -&+\\
- &  + & *&  +&- \\
+  &  - &  + & *&+\\
- &  + & -&  +& *  
\end{bmatrix},
\ee
where~$+$ denotes a non-negative entry,
$-$ denotes a non-positive entry, and~$*$ denotes ``don't care''. The matrix~$A^{[3]}$ is Metzler iff
\[
A=\begin{bmatrix}
*& + & 0& 0 & +  \\
+ & *& + & 0&0\\
0 &  + & *&  +&0 \\
0  &  0 &  + & *&+\\
+ &  0 & 0&  +& *  
\end{bmatrix},
\]
and~$A^{[2]}$ is Metzler iff
\[
A=\begin{bmatrix}
*& + & 0& 0 & -  \\
+ & *& + & 0&0\\
0 &  + & *&  +&0 \\
0  &  0 &  + & *&+\\
- &  0 & 0&  +& *  
\end{bmatrix}.
\]

\end{Example}

We consider the  three cases in Prop.~\ref{prop:sign_pattern_for_k_posi} in more detail.
In Case~\ref{case:nminus1},  $A^{[n-1]}$ is Metzler. Let~$U \in\{-1,1\}^{n\times n}$ 
be the antidiagonal matrix with entries  
\[
u_{ij}=\begin{cases}
(-1)^{j+1} , &\text{if } i+j=n+1,\\
0,& \text{ otherwise}.
\end{cases}
\]
Then~$U^T=U^{-1} $, and it can be shown that~$A^{[n-1]}$ Metzler implies that
that~$(-U  A U ^{-1} )$ is Metzler~\cite{Eyal_k_posi}. For example, when~$n=5$ a calculation gives 
\[
 U A U^{-1} = \begin{bmatrix}
                               a_{55} & -a_{54} & a_{53} &-a_{52} & a_{51} \\
                                - a_{45} & a_{44} &- a_{43} &a_{42} &- a_{41} \\ 
                                a_{35} & -a_{34} & a_{33} &-a_{32} & a_{31} \\
                               - a_{25} & a_{24} &- a_{23} &a_{22} &- a_{21} \\ 
                               a_{15} & -a_{14} & a_{13} &-a_{12} & a_{11} 
\end{bmatrix},
\] 
and~\eqref{eq:a45} implies that~$(-U  A U ^{-1} )$ is Metzler.
 In other words, 
in Case~\ref{case:nminus1}) 
the coordinate transformation~$y:=Ux$ gives~$\dot y= UAU^{-1}  y$, and this   is a competitive dynamical system~\cite{hlsmith}. 
Thus, $k$-positive systems, with~$k\in\{1,\dots,n-1\}$, may be viewed as a kind of interpolation from cooperative systems (when $k=1$) to competitive systems (when $k=n-1)$. 
  
\begin{figure}
 \begin{center}
  \includegraphics[scale=0.5]{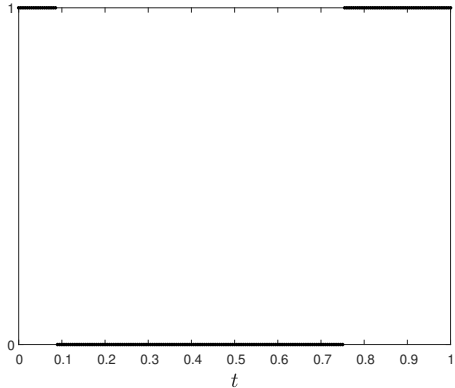}
	\caption{ $s^-(x(t,x_0))$ as a function of~$t$.}
	\label{fig:signs}
\end{center}
\end{figure}

In Case~\ref{case:kodd}),~$A$ is in particular Metzler. The dynamical behaviour in Case~\ref{case:keven})  is illustrated in the next example. 
 \begin{Example}
 Consider the case~$n=3$ and~$A=\begin{bmatrix}
 2  & 1  & -0.5\\ 
 0  & -1 &  0.5 \\ 
 -1 & 0  &   5 
 \end{bmatrix} $. Note that~$A$ is not Metzler, yet~$
 A^{[2]}=\begin{bmatrix}
 1   & 0.5 & 0.5 \\ 
 0   & 7   & 1   \\ 
 1   & 0   & 4
 \end{bmatrix}
$
 is Metzler (and irreducible). 
 Thm.~\ref{thm:k_pois_ltv} guarantees that for any~$x_0$ with~$s^-(x_0)\leq 1$, we have
 \be\label{eq:boundsmin}
 s^-(x(t,x_0))\leq 1\text{  for all }t\geq 0.
 \ee
 Fig.~\ref{fig:signs}
 depicts~$s^-(x(t,x_0))=s^-(\exp(At)x_0)$ for~$x_0=\begin{bmatrix}
 2 & -30 & -6
 \end{bmatrix}^T$. Note that~$s^-(x_0)=1$. It may be seen  that~$s^-(x(t,x_0))$ decreases and then increases, but always satisfies the bound~\eqref{eq:boundsmin}.
 \end{Example}
%
%
%
%

Fix~$ k\in\{2,\dots,n-1\} $. 
Prop.~\ref{prop:sign_pattern_for_k_posi}
implies that  a system 
  is~$k$-positive with~$k$ even 
iff it is~$2$-positive, and if it 
is~$k$-positive with~$k$ odd 
 then  it is~$1$-positive.

\subsection{Totally Positive Differential Systems}
We call~$A\in\R^{n\times n}$ a \emph{Jacobi matrix} if~$A$ is tri-diagonal with positive entries on the super- and sub-diagonals.
An immediate implication of Prop.~\ref{prop:sign_pattern_for_k_posi} is that~$A^{[k]}$ is Metzler and irreducible   for \emph{all}~$k\in\{1,\dots,n-1\}$ iff~$A$ is Jacobi. 
It then follows that for any~$t>0$ 
the matrices~$(\exp (At) )^{(k)} $, $k=1,\dots,n$, are (component-wise)  positive, that is,~$\exp(At)$ is TP for all~$t>0$. 
Combining this with Thm.~\ref{thm:k_pois_ltv} yields the following. 
\begin{Proposition}\cite{schwarz1970}\label{thm:TP}
The following two conditions are equivalent. 
\begin{enumerate} 
\item $A$ is Jacobi.
\item  for any~$x_0\in\R^n \setminus\{0\}$ 
the solution of the LTI~$\dot x(t)=Ax(t)$, $x(0)=x_0$,
satisfies
\[
s^-(x(t, x_0) ) \leq s^+(x(t, x_0) )\leq s^-(x_0) \text { for all } t>0.
\]
\end{enumerate}
\end{Proposition}

The left-hand side inequality follows from~\eqref{eq:sminussplus}, and the right-hand inequality   from 
the~SVDP of~TP systems in~\eqref{eq:SVDP}. Prop.~\ref{thm:TP} implies that~$s^-(x(t,x_0))$ 
and also~$s^+(x(t,x_0))$ are non-increasing  functions of~$t$, and may thus be considered as  piece-wise constant Lyapunov functions for the dynamics. 

Prop.~\ref{thm:TP} was proved by Schwarz~\cite{schwarz1970}, who only considered  linear systems. It was recently shown~\cite{margaliot2019revisiting} that important results on the asymptotic behaviour 
of time-invariant and periodic 
time-varying nonlinear systems with a Jacobian that is a Jacobi matrix for all~$t,x$~\cite{smillie,periodic_tridi_smith} follow from the fact that the associated  variational equation is a totally positive~LTV.

\section{$k$-Cooperative Systems}\label{sec:kcoop}
 
%
We now review the applications of~$k$-positivity
to the time-invariant nonlinear system:
\be\label{eq:time_invariant_non_linear}
\dot x=f(x),
\ee
with~$f\in C^1$. Let~$J(x):=\frac{\partial }{\partial x}f(x)$. We assume that the trajectories of~\eqref{eq:time_invariant_non_linear} evolve on a  convex and compact  state-space~$\Omega\subseteq\R^n$.

Recall that~\eqref{eq:time_invariant_non_linear}
is called \emph{cooperative} if~$J(x)$ is Metzler for all~$x\in \Omega$ \cite{hlsmith}. In other words, the variational equation associated with~\eqref{eq:time_invariant_non_linear}
is positive. The slightly stronger condition of strong cooperativity
has far reaching implications. By Hirsch's quasi-convergence theorem~\cite{hlsmith},
almost every bounded trajectory  
converges to the set of equilibria. 

It is natural to generalize cooperativity to~$k$-cooperativity by requiring  that the variational equation associated   
with~\eqref{eq:time_invariant_non_linear} is~$k$-positive.
 

\begin{Definition}\label{def:k-coop}\cite{Eyal_k_posi}
The nonlinear system~\eqref{eq:time_invariant_non_linear} is called  \emph{[strongly] $k$-cooperative} 
if the associated  LTV~\eqref{eq:var_eqn} is [strongly]~$k$-positive for any~$a,b\in\Omega$.
\end{Definition}

Note that for~$k=1$ this reduces to the definition of a cooperative  [strongly cooperative] dynamical system.

One immediate implication of Definition~\ref{def:k-coop} is the existence of certain invariant sets of the dynamics.

 
\begin{Proposition} \label{prop:inhgt}
Suppose that~\eqref{eq:time_invariant_non_linear}  
is~$k$-cooperative.
Pick~$a,b\in\Omega$. Then
\[
a-b\in P^k_- \implies  x(t,a)-x(t,b) \in P^k_-  \text{ for all } t \geq  0. 
\]
If, furthermore,~$0 \in \Omega $ and~$0$
is an equilibrium point of~\eqref{eq:time_invariant_non_linear}, i.e.~$f(0)=0$   then 
\[
a \in P^k_- \implies  x(t,a)  \in P^k_-  \text{ for all } t \geq  0. 
\]
\end{Proposition}

The sign pattern conditions in Prop.~\ref{prop:sign_pattern_for_k_posi}  yields  simple   sufficient conditions for
[strong] $k$-cooperativity  of~\eqref{eq:time_invariant_non_linear}. Indeed, if~$J(x)$ satisfies a sign pattern condition for all~$x\in \Omega$ then   the integral of~$J$ in the variational equation~\eqref{eq:var_eqn} satisfies the same sign pattern, and thus so does~$A^{ab}$. 
The next example, adapted from~\cite{Eyal_k_posi}, 
illustrates this.

\begin{Example} Elkhader~\cite{Elkhader1992}
studied the   nonlinear system
\begin{align}\label{eq:alexsys}
\dot x_1&=f_1(x_1,x_n),\nonumber\\
\dot x_i &= f_i(x_{i-1},x_i,x_{i+1}),
\quad i=2,\dots,n-1,\nonumber\\
\dot x_n&=f_n(x_{n-1},x_n),
\end{align}
under  the following assumptions:
the state-space~$\Omega\subseteq\R^n$
 is convex, $f_i\in C^{n-1}$, $i=1,\dots,n$, and 
there exist~$\delta_i\in\{-1,1\}$, $i=1,\dots,n$,
 such that
\begin{align*}
\delta_1\frac{\partial }{\partial x_n}f_1(x)   &>0,\\
\delta_2\frac{\partial }{\partial x_1}f_2(x) , \delta_3\frac{\partial }{\partial x_3}f_2(x)&>0,\\
&\vdots\\
\delta_{n-1} \frac{\partial }{\partial x_{n-2}}f_{n-1}(x),
\delta_n \frac{\partial }{\partial x_{n}}f_{n-1}(x)
 &>0,\\
\delta_n \frac{\partial }{\partial x_{n-1}}f_n(x)&>0,
\end{align*}
for all~$x\in\Omega$. 
This is a generalization of     the
monotone cyclic feedback system analyzed  in~\cite{poin_cyclic}.
As noted in~\cite{Elkhader1992}, we may assume without loss of generality that~$\delta_2=\delta_3=\dots=\delta_n=1$ and~$\delta_1 \in \{-1,1\}$. Then the Jacobian of~\eqref{eq:alexsys} has the form
\[
J(x)=\begin{bmatrix}
*& 0 &0 &0 &\dots & 0 & 0  & \sgn(\delta_1) \\
>0 & * &>0 &0 &\dots & 0& 0 & 0  \\
0& >0 & * &>0  &\dots &0&  0 & 0 \\
&&&\vdots\\
0&  0 & 0 & 0  &\dots & 0  &>0&  * \\
\end{bmatrix} ,
\]
for all~$x\in \Omega$. Here~$*$ denotes ``don't care''. Note that~$J(x)$ is irreducible for all~$x\in \Omega$. 

If~$\delta_1=1$ then~$J(x)$ is Metzler and Irreducible, so the system is strongly~$1$-cooperative. 
If~$\delta_1=-1$ then~$J(x)$ satisfies the sign pattern in Case~\ref{case:keven} in Prop.~\ref{prop:sign_pattern_for_k_posi}, so the system is strongly $2$-cooperative.
(If~$n$ is odd then~$J(x)$ also satisfies the sign pattern in Case~\ref{case:nminus1}, so  there is a coordinate transformation for which the  system is also strongly competitive.) 
\end{Example}

 The main result in~\cite{Eyal_k_posi} is that strongly~$2$-cooperative systems satisfy  
 a  
 strong Poincar\'{e}-Bendixson property.

 \begin{Theorem}  \label{thm:2dim}
	 (Poincar\'{e}-Bendixson Property)~\cite{9107214} Suppose that~\eqref{eq:time_invariant_non_linear} is strongly~$2$-cooperative. Pick~$a\in \Omega$.
If the omega limit set~$\omega(a)$ does not include an equilibrium 
then it is a closed orbit.
	\end{Theorem}

The proof of this  result is based on the seminal results of Sanchez~\cite{sanchez2009cones} on dynamical systems that admit an invariant cone of rank~$k$. Yet, Thm.~\ref{thm:2dim} is  considerably
stronger than  the main result in~\cite{sanchez2009cones}, as it applies to \emph{any} trajectory emanating from~$\Omega$   and not only to   \emph{pseudo-ordered} trajectories (see the definition in~\cite{sanchez2009cones}). 

The Poincar\'{e}-Bendixon property is useful because often it can be combined with a local analysis near  the equilibrium points to provide a global picture  of the dynamics. For a recent application of Thm.~\ref{thm:2dim} to   models from 
  systems biology and epidemiology,  see~\cite{margaliot2019compact,kordercont}.
  
Summarizing, the use of $k$-compounds allows to generalize positive linear systems [cooperative nonlinear systems] to $k$-positive linear systems [$k$-cooperative nonlinear systems]. In particular, $2$-cooperative systems, of any order $n$, ``behave like'' $2$-dimensional  systems.

The next section describes  the use of compound matrices to extend
the important concept of diagonal stability to $k$-diagonal stability~\cite{cheng_diag_stab}. We consider the case of discrete-time~(DT) systems.
For a symmetric matrix~$S\in\R^{n\times n}$, we use~$S\succ 0$ [$S\prec 0$]
to denote that~$S$ is positive-definite [negative-definite]. 
\section{$k$-Diagonal Stability}\label{sec:kdiag}
%
Recall that the discrete-time LTI system
\be\label{eq:dtlti}
x(j+1)=Ax(j),
\ee
is called \emph{diagonally stable} if there exists a diagonal matrix~$D \succ 0$ such that~$A^TDA\prec D$. 
In other words,~$V(x):=x^T D x$ is a
  diagonal Lyapunov function for~\eqref{eq:dtlti}.  
  
Diagonal stability of both discrete-time and continuous-time   LTI systems has attracted considerable interest, as it has important implications to non-linear dynamical systems~\cite{kaszkurewicz2012matrix,arcak2006diagonal,wimmer2009diagonal,barker1978positive}. For example, let~$S^1$
denote the class of functions~$f:\R^n \to\R^n$ such that~$f(x)=\begin{bmatrix}
f_1(x_1)&\dots&f_n(x_n))
\end{bmatrix}^T$,   each~$f_i:\R\to\R$ is continuous and satisfies
\[
0<|f_i(s)|\leq |s|, \text{ for all } s\in\R\setminus\{0\}.
\]
Note that~$f(0)=0$. Consider the DT non-linear system
\be\label{eq:dtnonline}
x(j+1)=A f(x(j)),\quad f\in S^1 .
\ee
If~\eqref{eq:dtlti} is diagonally stable then it is straightforward to verify that~$V(x):=x^T D x$ is also a Lyapunov function for the nonlinear system~\eqref{eq:dtnonline}.

 Stable LTIs always admit a quadratic Lyapunov function, but not necessarily a diagonal Lyapunov function~(DLF). 
However, stable \emph{positive} LTIs always  admit a DLF. 
\begin{Lemma}\label{lemma:diagonal_stability_pos_DT_LTI}\cite{berman87,  RANTZER201572} 
%
If $A\in\mathbb{R}^{n\times n}$ with $A\ge 0$
then the following statements are equivalent:
\begin{enumerate}
    \item The matrix $A$ is Schur, i.e., $\rho (A) < 1$.
    \item \label{item:xi} There exists $\xi \in\mathbb{R}^{n}$ with~$\xi\gg 0$
    such that $A\xi \ll \xi$.
    \item  \label{item:zed} There exists $z\in\mathbb{R}^{n}$ with~$z\gg 0$
    such that $A^{T} z \ll z$.
    \item \label{item:diagonald}
    There exists a diagonal 
    matrix~$D \succ  0  $
    such that~$A^{T} DA \prec D$.
    \item The matrix~$(I - A)$ is non-singular and~$(I - A)^{-1}  \geq 0$.
\end{enumerate}
%
\end{Lemma}
%
%
\begin{Remark}
%
Suppose that~$A\geq 0$ is Schur. 
Fix~$x,y\in\mathbb{R}^{n}$,
with~$x,y \gg 0$. Then
\begin{align*}
\xi& := (I-A)^{-1} x,\\
  z& := (I-A^{T})^{-1} y ,\\
 D& := \diag \left ( \frac{z_1}{\xi_1},\hdots,
\frac{z_n}{\xi_n} \right )
\end{align*}
satisfy conditions~\ref{item:xi}),
\ref{item:zed}), and~\ref{item:diagonald})
in 
Lemma~\ref{lemma:diagonal_stability_pos_DT_LTI},
respectively. This provides a constructive recipe for determining the DLF~$D$. 
Note that if~$A\in\mathbb{R}^{n\times n}$ 
is Schur and~$A\le 0$, then $(-A)$ is
a non-negative Schur matrix, and 
Lemma~\ref{lemma:diagonal_stability_pos_DT_LTI} guarantees   that there exists a diagonal 
matrix~$D\succ 0 $   such that $A^{T} DA \prec D$.
%
\end{Remark}

\subsection{Discrete-Time~$k$-Diagonal Stability}
Fix~$k\in\{1,\dots,n\}$.
To study  the evolution     of volumes of~$k$-parallelotopes under the DT LTI~\eqref{eq:dtlti}, fix~$k$
initial conditions~$a^1,\dots,a^k \in \R^n$. Then
\begin{align}\label{eq:comp_dos}
    \begin{bmatrix}
    x(j,a^1)&\dots&x(j,a^k)
    \end{bmatrix}^{(k)}& =  \begin{bmatrix}
    A^j a^1&\dots&A^j a^k
    \end{bmatrix}^{(k)}\nonumber \\
    &= (A^j)^{(k)} \begin{bmatrix}
      a^1&\dots&  a^k
    \end{bmatrix}^{(k)}\nonumber \\
    &=(A^{(k)})^j \begin{bmatrix}
      a^1&\dots&  a^k
    \end{bmatrix}^{(k)}.
\end{align}
This may be interpreted as the~$k$-compound system of~\eqref{eq:dtlti}.
 It shows that  the evolution of~$k$-parallelotopes  
follows  a  DT LTI with the matrix~$A^{(k)}$.  This naturally leads to the following definition. 
\begin{Definition}\cite{cheng_diag_stab}
Consider
the DT LTI~\eqref{eq:dtlti} with~$A\in\R^{n\times n}$. Fix~$k\in\{1,\dots,n\}$ and let~$r:=\binom{n}{k}$. Then~\eqref{eq:dtlti}
is called \emph{$k$-diagonally stable} 
if there exists a   matrix~$D =\diag(d_1,\dots,d_r)\succ 0 $
such that
\begin{align}\label{eq:diagonal_stability}
(A^{(k)})^T DA^{(k)}  \prec D. 
\end{align}
\end{Definition}
Note that for~$k=1$ this  reduces to   standard   diagonal stability, and for~$k=n$ this becomes   the requirement~$(\det(A))^2<1$. 
Note also that if~$A^{(k)} $ is Schur and~$A^{(k)}\geq 0 $ 
then Lemma~\ref{lemma:diagonal_stability_pos_DT_LTI} implies that~\eqref{eq:dtlti}
is $k$-diagonally stable.



It is natural to expect that diagonal stability implies~$k$ diagonal stability for any~$k$. The next result shows that this is indeed so. 
  \begin{Proposition}\cite{cheng_diag_stab}
  If~$A $ is diagonally stable then it is $k$-diagonally stable for any~$k\in\{1,2,\dots,n\}$. 
  \end{Proposition}
\begin{IEEEproof}
%
Let~$\mathbb{D}^{n\times n}$ denote  the set of positive diagonal~$n\times n$ matrices. 
There exists a~$P\in\mathbb{D}^{n\times n}$
such that~$A^{T} P A \prec P$.
Multiplying this inequality by~$P^{-1/2}$ on the right and left-hand side     yields:
\be\label{eq:pm1/2}
P^{-1/2} A^{T} P A P ^{-1/2} \prec I_{n},
\ee
so  the symmetric matrix~$
  H:=  P^{-1/2} A^{T} P A P ^{-1/2}
$
is Schur.
Fix~$k\in \{1,2,\hdots,n\}$, and let~$D:= P^{(k)}$.
Then~$D\in \mathbb{D}^{r\times r}$, where~$r:= \binom{n}{k}$. Now, 
\begin{align*}
     H ^{(k)}
    =
    D^{-1/2} ( A^{(k)} )^{T} D A^{(k)} D^{-1/2},
\end{align*}
and since~$ H ^{(k)}$ is symmetric and any eigenvalue of~$ H ^{(k)}$
is the product of~$k$ eigenvalues of~$H$, we also have~$ H ^{(k)}
    \prec I_{r}$. We conclude that 
\begin{align*}
    (A^{(k)})^{T} D A^{(k)} \prec D,
\end{align*}
i.e.,~$A $ is $k$-diagonally stable.
\end{IEEEproof}

We note that, in general, $k$-diagonal stability, with~$k>1$, does not imply diagonal stability (see the  specific examples in~\cite{cheng_diag_stab}).


It is interesting to find general classes of matrices  that are $k$-diagonally stable for some~$k$.  
Recall that~$A\in\mathbb{R}^{n\times n}$
is called cyclic if 
\begin{align*}
    A = 
    \begin{bmatrix}
    \alpha_1 & \beta_1 & 0 & \hdots & 0& 0 \\
    0 & \alpha_2 & \beta_2 & \hdots & 0 & 0 \\
    0 & 0 & \alpha_3 & \hdots & 0 & 0 \\
    \vdots & \vdots & \vdots & \ddots & \vdots& \vdots \\
    0 & 0 & 0 & \hdots &\alpha_{n-1}& \beta_{n-1} \\
    (-1)^{\ell +1}\beta_{n} & 0 & 0 & \hdots& 0 & \alpha_{n}
    \end{bmatrix}
\end{align*}
with $\alpha_{i},\beta_{i} \ge 0$, $i=1,\hdots,n$,
and $\ell \in\mathbb{Z}$.
%
%
The DT LTI system $x (j+1) = A x(j)$
is called cyclic if~$A$ is cyclic.
This represents a dynamical system where each~$x_{i}$, $i=1,\dots, n-1$,
receives  positive feedback  from its ``cyclic neigbours''~$x_{i} , x_{i+1} $,
with the exception that~$x_{n} $ receives  negative  feedback 
from~$x_{1} $
  if $\ell$ 
is even. 

\begin{Theorem}\label{Thm:Cyclic_Theorem}\cite{cheng_diag_stab}
%
Suppose that~$A\in\R^{n\times n}$  
is cyclic for some
$\ell\in \{1,\hdots,n-1\}$.
Then~$A$ is~$SR_{\ell}$
with signature~$\epsilon_{\ell} = 1$.
Furthermore, 
if $\ell$ is odd, 
then~$A$ is   diagonally 
stable 
iff $A$ is Schur. If $\ell$ is even, 
then~$A$ is   
$\ell$-diagonally stable 
iff~$A^{ (\ell) }$ is Schur.
%
\end{Theorem}
%
\begin{Example}
Consider the case~$n=3$, i.e., 
\begin{align*}
    A = 
    \begin{bmatrix}
    \alpha_1 & \beta_1 & 0 \\
    0 & \alpha_2 & \beta_2 \\ 
    (-1)^{ \ell + 1 } \beta_3 & 0 & \alpha_3
    \end{bmatrix}.
\end{align*}
Then
\begin{align*}
    A^{(2)} =
    \begin{bmatrix}
    \alpha_1 \alpha_2 & \alpha_1 \beta_2 & \beta_1 \beta_2 \\
    (-1)^{\ell} \beta_1 \beta_3 & \alpha_1 \alpha_3 & \alpha_3 \beta_1 \\
    (-1)^{\ell} \alpha_2 \beta_3 & (-1)^{\ell} \beta_2 \beta_3 & \alpha_2 \alpha_3
    \end{bmatrix}.
\end{align*}
Thus, for~$\ell=1\, [\ell=2]$,~$A$ is~$SR_1$ $[SR_2]$ with signature~$\epsilon_1=1$ $[\epsilon_2=1]$ since all the entries of~$A$ $[A^{(2)}]$ are non-negative.
\end{Example}

As an application of these notions,  Ref.~\cite{cheng_diag_stab} describes a class  of non-linear functions~$S^k$ such that~$k$-diagonal stability of~\eqref{eq:dtlti}
implies that   the evolution of~$k$-parallelotopes under the non-linear dynamics
\[
x(j+1)=Af(x(j)), \quad f\in S^k,
\]
is asymptotically stable. 
For~$k=1$, $S^k$ is the set of functions~$S^1$ defined above.

The $k$-compound matrices are   defined only for integer values of~$k$, as they are based on~$k\times k$ sub-matrices. However, it turns out that there are good reasons to   generalize the notion of~$k$-compounds, with~$k$ an integer, to~$\alpha$-compounds, with~$\alpha$ a positive \emph{real} number. In particular, this allows to introduce the notion of~$\alpha$-contracting systems. 
%
\section{$\alpha$-Compounds and $\alpha$-Contracting~Systems}\label{sec:alpha}
%
A  recent paper~\cite{wu2020generalization} defined a generalization
called the~$\alpha$-multiplicative compound and~$\alpha$-additive compound of a matrix, where~$\alpha$ is a  {real} number.

Recall that the  Kronecker product
 of 
two matrices~$A \in \C^{n \times m}$ and~$B \in \C^{p \times q}$ is
\be
A \otimes B := 
\begin{bmatrix}
a_{11} B & a_{12}B & \cdots & a_{1m} B  \\
a_{21} B & a_{22}B & \cdots & a_{2m} B \\
\vdots & \vdots  & \ddots & \vdots \\
a_{n1} B & a_{n2}B & \cdots& a_{nm} B
\end{bmatrix},
\ee
where $a_{ij}$ denotes the $ij$th entry of $A$. Hence, $A \otimes B \in \C^{(np) \times (mq)}$. 
The Kronecker sum of two
square matrices~$X \in \C^{n \times n}$ and~$Y \in \C^{m \times m}$ is
 \begin{equation}
    X \oplus Y :=  X \otimes I_m + I_n \otimes Y.
 \end{equation}

\begin{Definition}
Let~$A\in\C^{n\times n} $ be non-singular. 
If~$\alpha=k+s$, where~$k\in\{1,\dots,n-1\}$ and~$s\in(0,1)$ then the~$\alpha$-multiplicative compound of~$A$ is defined by:
\be\label{Eq:defaalf}
A^{(\alpha)} : = (A^{(k)})^{1-s} \otimes (A^{(k+1)})^{s}.
\ee
\end{Definition}

Note that $A^{(\alpha)} \in \C^{r \times r}$, where~$r : = \binom{n}{k}\binom{n}{k+1}$, and that~$A^{(\alpha)}$ may be complex (non-real) even if~$A$ is real. Since $A$  is non-singular,    $A^{(\ell)}$ is   non-singular for all~$\ell \in \{1, \dots, n\}$, so~$(A^{(k)})^{1-s}$ and~$(A^{(k+1)})^s$ in~\eqref{Eq:defaalf} are well-defined.

 The matrix~$A^{(\alpha)}$  is a kind of ``multiplicative interpolation'' between~$ A^{(k)} $ and~$A^{(k+1)} $.
For example,~$A^{(2.2)}    = (A^{(2)})^{0.8} \otimes (A^{(3)})^{0.2}$.

\begin{Example}\label{exa:diagexa}
Let~$D=\diag(d_1,\dots,d_4)$ with~$d_i\not=0$ for all~$i$. Fix~$\alpha  \in (2,3)$, so that~$k=2$ and~$s=\alpha-2 \in (0,1)$. 
Then
\begin{align*}
    D^{(\alpha)}&=( D ^{(2)})^{1-s} \otimes (D^{(3)})^s\\&=
   \diag( (d_1 d_2)^{1-s},(d_1 d_3)^{1-s},\dots,
    (d_3 d_4)^{1-s} )
   )\\& \;\;
   \otimes
   \diag(
      (d_1 d_2 d_3)^{s},
      (d_1 d_2 d_4)^{s}, 
      (d_1 d_3 d_4)^{s}, 
      (d_2 d_3 d_4)^{s}
   )
   \\&=
   \diag( d_1 d_2 d_3^s, d_1 d_2 d_4^s, \dots,    d_2^s d_3 d_4 ) , 
\end{align*}
so, any eigenvalue of~$D^{(\alpha)}$ is a  ``multiplicative interpolation'' between   eigenvalues of~$ D ^{(2)}$ and~$ D ^{(3)}$.
\end{Example}

Just like the~$k$-additive compound, 
the~$\alpha$-additive compound is defined
using the~$\alpha$-multiplicative compound.

\begin{Definition}
Let~$A\in\C^{n\times n} $ be non-singular. 
If~$\alpha=k+s$, where~$k\in\{1,\dots,n-1\}$ and~$s\in(0,1)$ 
then the~$\alpha$~additive compound matrix of~$A $ is 
\begin{align*}
    A^{[\alpha]} := \frac{d}{d\epsilon} (I+\epsilon A)^{(\alpha)} |_{\epsilon=0} , 
\end{align*}
\end{Definition}
It was shown in~\cite{wu2020generalization}
that this yields
\be\label{eq:aalpha}
 A^{[\alpha]} = ((1-s)A^{[k]}) \oplus (sA^{[k+1]}).
 \ee

\begin{Example}\label{eq:aalpnm1}
Let~$A\in \R^{n\times n}$ and fix~$\alpha \in (n-1,n)$ so that~$\alpha=k+s$, with~$k=n-1$ and~$s\in(0,1)$. 
Then~\eqref{eq:aalpha} gives
\begin{align}\label{eq:anmqk}
 A^{[\alpha]} &= ((1-s)A^{[n-1]}) \oplus (sA^{[n]}) \nonumber\\
&=((1-s)A^{[n-1]}) \oplus (s \trace(A))\nonumber\\
&=((1-s)A^{[n-1]}) \otimes I_ 1 +I_n \otimes  (s \trace(A))\nonumber\\
&=(1-s)A^{[n-1]}  +   s \trace(A)  I_n ,
\end{align}
so in this particular case the Kronecker sum
simplifies to a standard 
matrix sum. 
\end{Example}

  Consider the time-varying non-linear  dynamical system~\eqref{eq:nonlin} where~$f$ is~$C^1$.
  Let~$x(t,t_0,x_0)$ denote the solution of~\eqref{eq:nonlin}  at time~$t$ with~$x(t_0)=x_0$. We assume from here on that~$t_0=0$, and let~$x(t,x_0):=x(t,0,x_0)$. We also assume that  the system admits  a convex  invariant set~$\Omega \subseteq \R^n$, that is, for any~$x_0\in \Omega$ we have~$x(t,x_0) \in \Omega $ for  all~$t\geq 0$. Let~$
  J (t,x):=\frac{\partial }{\partial x}f(t,x),
 $
   
Recall that the $k$-additive compound was used to define the notion of~$k$-contraction. The next definition follows  the same line of reasoning. 

\begin{Definition} \cite{wu2020generalization}
The system~\eqref{eq:nonlin}
is  called \emph{infinitesimally  $\alpha$-contracting} w.r.t. the norm~$|\cdot| $   if
\be \label{eq:alphacon}
\mu(J^{[\alpha]}(t,x))\leq -\eta<0 ,
\ee
for all~$t\geq 0$ and~all~$x$ in the state space~\cite{wu2020generalization}. 
\end{Definition}

 A set~$K\subseteq \Omega$ is called a strongly invariant set of~\eqref{eq:nonlin}
  if
  \be\label{eq:invtau} 
  K =  x(t, K) \text{ for all } t\geq 0.
  \ee
  In other words, if we take 
  all the  points in~$K$ 
  as initial conditions of the dynamics then the  union  of
  the resulting solutions at time~$t$ is~$K$. 
  For example, an equilibrium or a limit cycle are strongly  invariant  sets. More generally, for any~$a\in \Omega$, the omega limit set~$\omega(a)$ is a   strongly invariant set.

Using these notions, it is possible to restate a  seminal result  of  Douady 
 and Oesterl\'{e}~\cite{Douady1980} (see also~\cite{smith_hauss})  in terms of~$\alpha$
 contraction. 
\begin{Theorem}\label{thm:alpha}\cite{wu2020generalization}
Suppose that~\eqref{eq:nonlin}
is~$\alpha$-contracting for some~$\alpha\in [1,n]$.
Then any compact and strongly invariant set of the dynamics 
has a  Hausdorff dimension smaller than~$\alpha$. 
\end{Theorem}

Roughly speaking, the dynamics contracts   sets with
a  larger Hausdorff dimension,
so any strongly invariant set must have a Hausdorff dimension small than $\alpha$.

\begin{figure}[t]
 \begin{center}
\includegraphics[width=8cm,height=6cm]{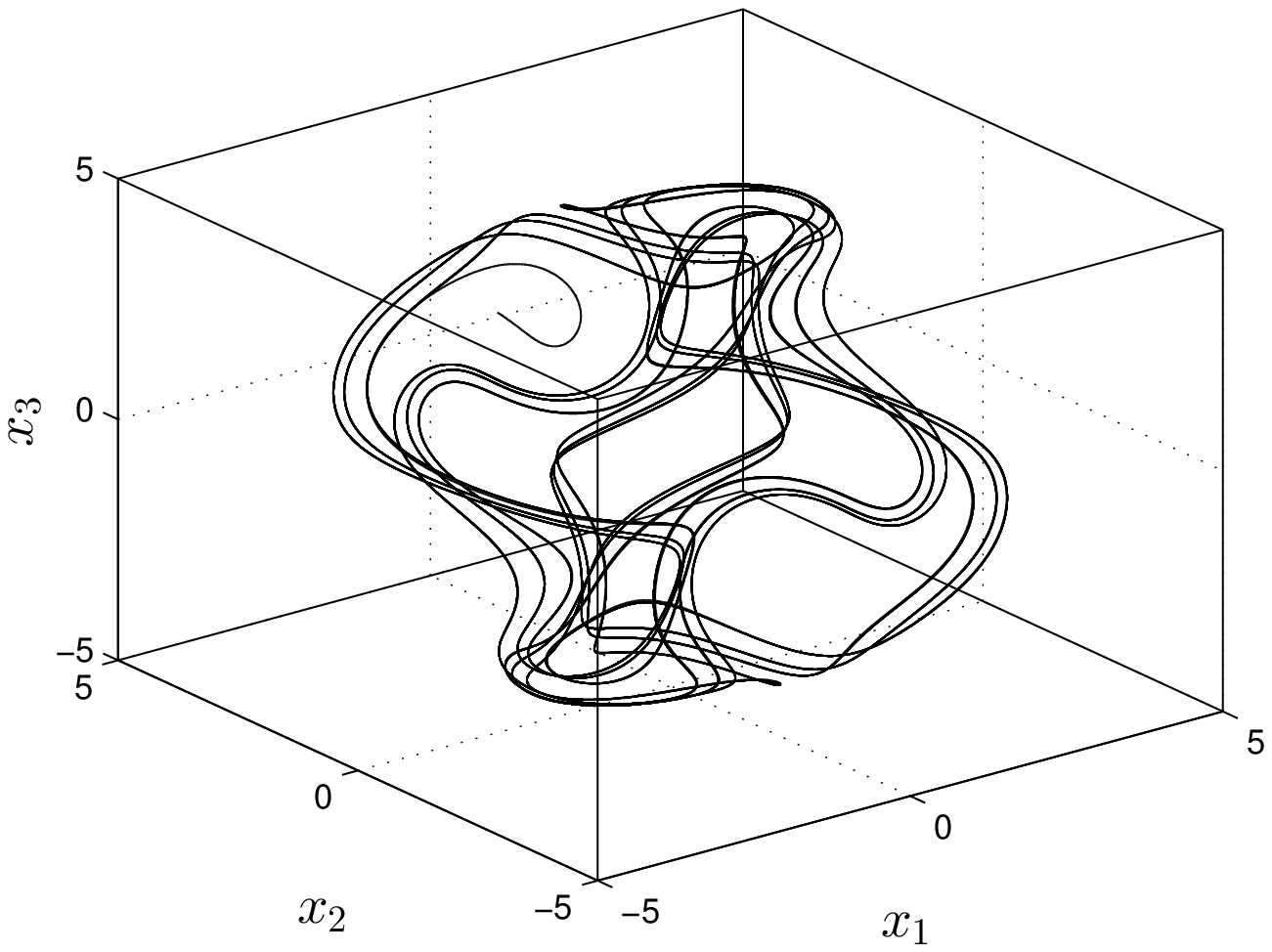}
\caption{A trajectory of~\eqref{eq:thom} with~$b=0.1$  emanating from~$x(0)=\begin{bmatrix} 
  1 & -2 & 1 \end{bmatrix}^T$.
}\label{fig:chaos}
\end{center}
\end{figure}

\begin{figure}[t]
 \begin{center}
\includegraphics[width=8cm,height=6cm]{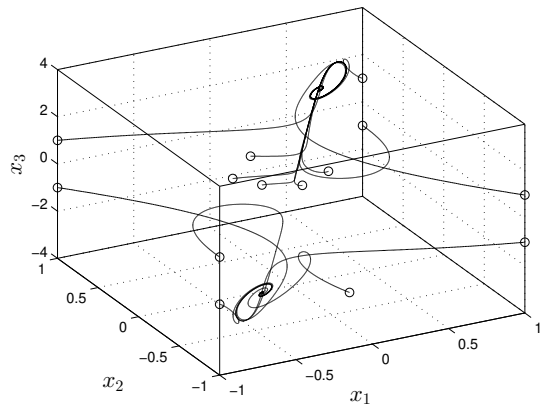}
\caption{Several trajectories of  the closed-loop system.  The circles denote the initial conditions of the trajectories.  }\label{fig:chaos_closed}
\end{center}
\end{figure}

The next example, adapted from~\cite{wu2020generalization}, shows how these notions  can be used to design a feedback controller that 
``de-chaotifies'' a nonlinear dynamical system. 
\begin{Example}
Thomas' cyclically symmetric attractor~\cite{thomas99,chaos_survey} is
a popular example for a chaotic system. It is described by:
\begin{align} \label{eq:thom}
\dot x_1 =&  \sin(x_2)-bx_1 , \nonumber \\
\dot x_2 =&  \sin(x_3) - bx_2, \\
\dot x_3 =& \sin(x_1) - bx_3, \nonumber
\end{align}
where $b>0$ is the
dissipation  constant. 
The convex and compact set~$\Omega : = \{x\in\R^3:  \displaystyle \max_{i }|x_i|
  \leq b^{-1} \}$ is an
invariant set of the dynamics, and we consider only initial conditions in this set.

For~$b > 1$ the origin is the single stable equilibrium  of~\eqref{eq:thom}. As~$b$ is decreased, the dynamics becomes more complicated. 
Fig.~\ref{fig:chaos}
 depicts  the solution of the system emanating from~$x(0) = \begin{bmatrix} 1& -2&1\end{bmatrix}^T$ for~$ b=0.1$.
 Note the symmetric strange attractor. 
 
The Jacobian $J_f$ of the vector field in~\eqref{eq:thom} is
\begin{align*}
   J_f (x)=\begin{bmatrix}-b&\cos(x_2)&0 \\ 0&-b&\cos(x_3) \\ \cos(x_1)&0&-b\end{bmatrix},
   \end{align*}
   and thus~$
   J_f ^{[3]}=\trace(J (x))=-3b$. 
Since~$b>0$, this implies that the system is $3$-contracting  w.r.t. any norm. Let~$\alpha = 2+s$, with~$s\in(0,1)$. Then combining~\eqref{eq:anmqk} 
and~\eqref{eq:add23}
gives
\begin{align*}
    J_f ^{[\alpha]}(x)
    &=(1-s)J_f^{[2]}(x)\oplus sJ_f^{[3]}(x)\\
    &=
    \begin{bmatrix}
    -(2+s)b & (1-s)\cos(x_3) & 0 \\
    0 & -(2+s)b & (1-s)\cos(x_2) \\
    -(1-s)\cos(x_1) & 0 & -(2+s)b \\
    \end{bmatrix}.
\end{align*}
Thus,
\be\nonumber
\mu_{1 }(J_f ^{[\alpha]}(x)) \leq 1-2b-s(b+1), \text{ for all } x\in \Omega .
\ee
We conclude that for any~$b\in(0,1/2) $
 the system is~$(2+s)$-contracting w.r.t. the~$L_1$ norm for any
$
    s>\frac{1-2b}{1+b}
$.

We now show how $\alpha$-contraction can be used to design a partial-state   controller for the system guaranteeing that the closed-loop system has a ``well-ordered'' behaviour. 
Suppose that 
the closed-loop system is:
\begin{align}\label{eq:cldloopsystem}
\dot x = f(x) + g(x),
\end{align}
where~$f$ is the vector field in~$\eqref{eq:thom}$ and~$g$ is the controller.
Let~$\alpha = 2 + s$, with~$s \in(0,1)$. 
The Jacobian of the closed-loop system is~$J_{cl}:=J_f+J_g$,  so
\begin{align*}
    \mu_1(J_{cl}^{[\alpha]})& = \mu_1(J_f^{[\alpha]}+J_g^{[\alpha]}) \\& \leq \mu_1(J_f^{[\alpha]})+\mu_1(J_g^{[\alpha]}) \\
    &\leq 1-2b-s(b+1)+\mu_1(J_g^{[\alpha]}),
\end{align*}
where the first  inequality follows from the fact that any matrix measure is sub-additive. This implies that the closed-loop system is~$\alpha$-contracting if
\begin{align}\label{eq:cond_cont}
    \mu_1(J^{[\alpha]}_g (x) ) <    s(b+1)+2b-1
    \text{ for all } x\in \Omega  . 
\end{align}
Consider, for example,   the controller~$g(x_1,x_2):=c \diag(1,1,0) x $, with gain~$c<0$.
Then~$J_g^{[\alpha]} = c\diag(  2  ,1+ s,1+ s )$
and for any~$c<0$ condition~\eqref{eq:cond_cont} becomes
\begin{align}\label{eq:cond_cont1}
     (1+ s)c <    s(b+1)+2b-1.
\end{align}
This provides a simple recipe  for 
determining the gain~$c$ so that the closed-loop system is~$(2+s)$-contracting. For example,  
  when~$s \to 0$, Eq.~\eqref{eq:cond_cont1} yields
$
c<2b-1
$,
and this guarantees that the closed-loop system is~$2$-contracting. Recall that in a time-invariant~$2$-contracting system 
every bounded 
trajectory converges to the set of equilibria,
 thus ruling out chaotic attractors and even non-trivial limit cycles~\cite{li1995}.  
Fig.~\ref{fig:chaos_closed} depicts the behaviour of the closed-loop system \eqref{eq:cldloopsystem} with~$b=0.15$   and~$c=2b-1.15$. The closed-loop system is thus $2$-contracting,
and as expected  
every solution converges to an equilibrium. 

Summarizing, the notion of   $\alpha$-contraction allows to add to the chaotic system the ``correct amount'' of feedback contraction, that may be a \emph{non-integer}, to obtain a~$2$-contracting system. 
\end{Example}

\begin{figure}
 \begin{center}
  \includegraphics[scale=0.8]{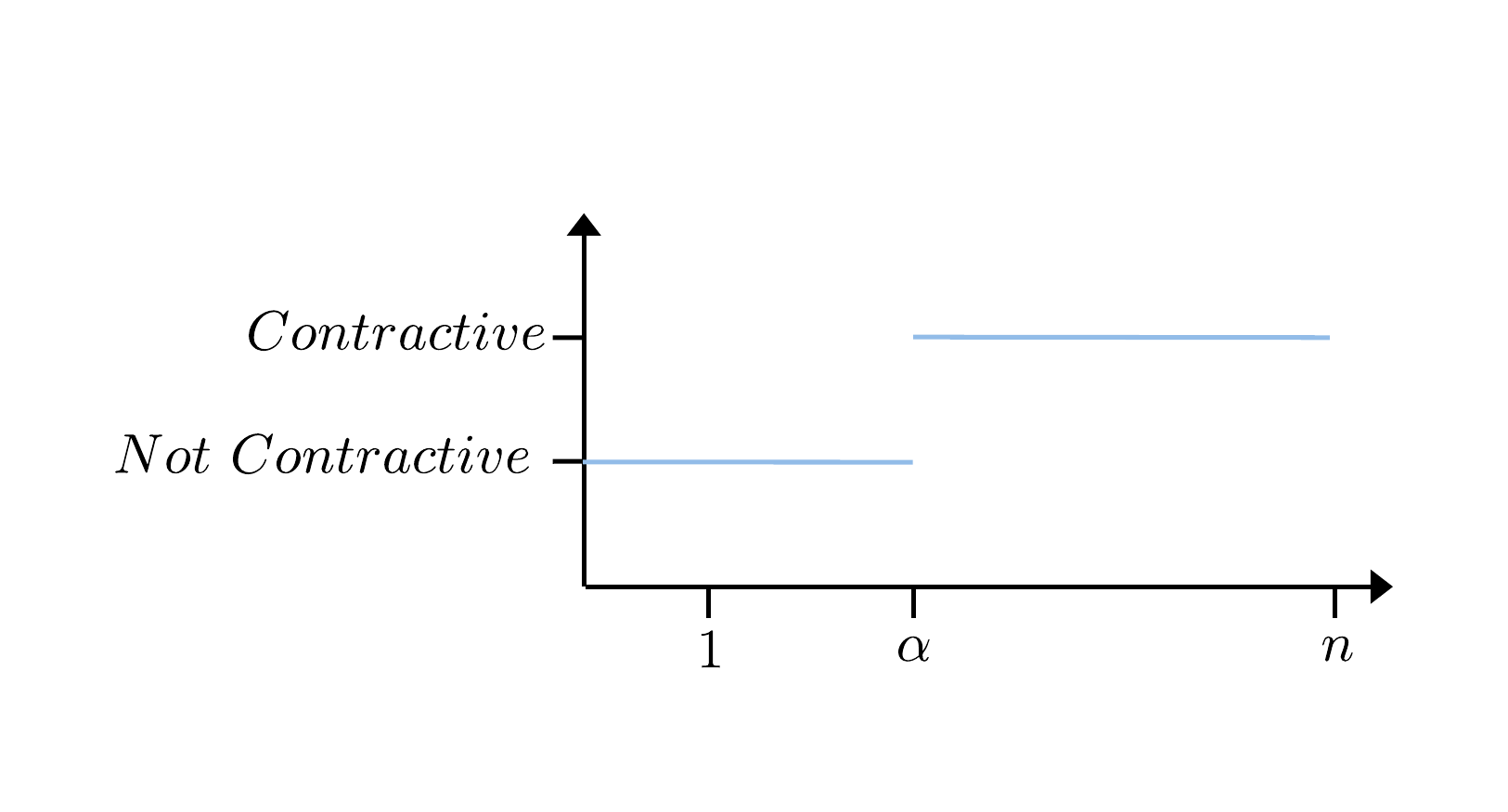}
  \caption{$\alpha$-contraction is  not a discrete property like~$k$-contraction.   }
    \label{fig:alpha_contr}
\end{center}
\end{figure}

The notion of
$\alpha$-contraction has another interesting implication. It was shown in~\cite{kordercont} that if
\[
\mu(J^{[\alpha]})\leq-\eta<0,
\]
with~$\mu$ a  matrix measure induced by an~$L_p$ norm, with~$p\in\{1,2,\infty\}$, then
\[
\mu(J^{[\alpha+\varepsilon]})\leq-\eta<0, \text{ for any } \varepsilon\geq 0.
\] 
This monotonicity 
property implies 
the following. If the system is~$(n-1)$-contracting then 
there always exists a \emph{real} value~$\alpha^*$ such that the system is~$\alpha $ contracting for any~$\alpha>\alpha^*$. 
 This implies that contraction is not a binary, yes-no property, that is, a system is either contracting  or not. 
Rather, contraction holds for any~$\alpha >\alpha^*$ (see Fig.~\ref{fig:alpha_contr}).

The next example demonstrates this in the context of the well-known 
linear consensus algorithm~\cite{eger_consens}. 
\begin{Example}
Consider the LTI system
\be\label{eq:lap}
\dot x= -Lx,
\ee
  where~$L\in\R^{n\times n}$ is   the Laplacian
 of a   (directed or undirected) weighted graph with a globally reachable vertex.
For any~$c\in\R$ we have that~$c 1_n$ is an equilibrium of~\eqref{eq:lap}, so the system cannot be $1$-contracting w.r.t. any norm. 

We now show that  for \emph{any}~$s >0$  
  there exists a  vector norm~$|\cdot|$ such that~\eqref{eq:lap}
is~$(1+s)$-contracting w.r.t. to~$|\cdot|$ (in other words, 
although the system is not~$1$-contracting  it is a ``arbitrarily close'' to being~$1$-contracting). 
It is enough to show that this holds for any~$s\in(0,1)$. Fix such an~$s$.
Let~$\lambda_i$ denote the eigenvalues of~$L$ ordered
such that
\begin{align} \label{eq:ord_eig}
\real(\lambda_1) \leq \real(\lambda_2) \leq \cdots \leq \real(\lambda_n), 
\end{align}
Then~$\lambda_1=0$
and~$\real(\lambda_2)>0$.
Eq.~\eqref{eq:aalpha}
gives  %
\begin{align*}
 L^{[1+s]} &= ((1-s) L ) \oplus (s L^{[2]}) . 
 \end{align*}
 The eigenvalues of~$L^{[2]}$ are the sums~$\lambda_i+\lambda_j$, $i<j$, so in particular the eigenvalue of~$L^{[2]} $
 with the minimal  real part is~$\lambda_2$.
Recall
that the eigenvalues 
of~$A\oplus B$ are
  the pairwise sums of the spectrum 
of~$A$ and~$B$,
so the minimal  real part 
of an eigenvalue of~$L^{[1+s]} $
  is~$s\real(\lambda_2)> 0$, i.e.~$-L^{[1+s]}$ is Hurwitz. It is well-known~\cite{sontag_cotraction_tutorial,Wang2add} that this implies that there exists a matrix measure~$\mu$, induced by a scaled~$L_2$ norm,
such that~$\mu (- L^{[1+s]})<0$.
Thus,~\eqref{eq:lap} is~$(1+s)$-contracting. 

Note that combining this with  Thm.~\ref{thm:alpha}
implies that 
 any compact and strongly invariant set of~\eqref{eq:lap}
has a  Hausdorff dimension smaller or equal to  one. Indeed, the  one-dimensional
``diagonal set''~$\{x\in\R^n\st x_1=\dots=x_n   \}$ is an   invariant set of the dynamics.
\end{Example}

So far we only considered autonomous systems  i.e., systems with no inputs. The next section briefly  reviews recent work by Grussler, Sepulchre  and their colleagues on 
 single-input single-output~(SISO) DT LTI systems
 that are $k$-positive. This leads to the interesting notion of the \emph{$k$-compound of a transfer function}. For the sake of consistency in this paper,  we modify some of the    terminology and notations.
Recall that a square matrix~$A$ is called
Hankel if every entry~$a_{ij}$ depends only on~$i+j$. 
In other words, the entries along any anti-diagonal are equal. 

\section{Hankel $k$-Positivity  of Discrete-Time SISO LTI Systems}  \label{sec:hankel}
%
Consider the  SISO  DT LTI system
\begin{align}\label{eq:dtisos}
     x(j+1) &= Ax(j) + bu(j), \nonumber \\
    y(j) &= c^{T} x(j),
\end{align}
with~$A\in\R^{n\times n}$, and~$b,c\in\R^n$.
The impulse response of this system      is:
\begin{align}\label{eq:impulsegj}
    g(j): = c^{T} A^{j-1} b\; \mathbbm{1}(j-1),
\end{align}
where $\mathbbm{1}(\cdot)$ is the Heaviside  step function.
We assume throughout that~$\sum_k |g(k) | <\infty$, and consider controls
satisfying~$\sup_k |u(k)| <  \infty$.

Let
\begin{align*}
G(z) := \sum_{k=0}^{\infty} g(k) z^{-k}
\end{align*}
denote the transfer function of~\eqref{eq:dtisos}. We always assume that the state-space representation~\eqref{eq:dtisos}
is minimal, so there are no zero-pole cancellations in~$G(z)$.

Section~\ref{sec:posi} considered systems without inputs
satisfying a  SVDP w.r.t. to the mapping from a past state vector to a future state vector.   
For an input-output system 
 like~\eqref{eq:dtisos}  it is natural to consider the relation between the number of sign variations
 in a sequence of inputs
and the number of sign variations 
in the  corresponding sequence of outputs. In certain applications  e.g.,   smoothing filters, it is natural to require an SVDP between these sequences.
In general, establishing such an SVDP  for~\eqref{eq:dtisos}
is a difficult  problem because any sequence of outputs is the sum of two terms: the first [second] representing  the effect of the initial condition [input sequence]. 
Bounding the number of sign changes in the sum of two vectors is non-trivial. For example, 
for~$a =\begin{bmatrix}
2&2&2 
\end{bmatrix}^T$
and
$
b =\begin{bmatrix}
-1&-3&-1& 
\end{bmatrix}^T$, 
we have~$s^-(a)=s^-(b)=0$, yet~$s^-(a+b)=  2 $.

Thus,  it is customary to make some simplifying assumptions. From here on, we consider   the so-called \emph{Hankel case} defined by the assumptions that
\be\label{eq:uposi}
0=u(0)=u(1)=u(2)=\dots,
\ee
and that the relevant sub-sequence of outputs  is
\[
y(0),\;y(1),\;y(2),\dots. 
\]
Intuitively speaking, the control values~$u(\ell)$, with~$\ell<0$, 
determine an initial condition~$x(0)$, and~\eqref{eq:uposi} implies that 
for any~$j\geq 0$ the system~\eqref{eq:dtisos} reduces to the autonomous system
\begin{align*} 
     x(j+1) &= Ax(j),   \\
    y(j) &= c^{T} x(j).
\end{align*}
%

%
%
%
%
%
Under these assumptions, for any~$j\geq 0$ we have 
\begin{align}\label{eq:conv}
    y(j)   =
    ( H_{g} u ) (j) :=
    \sum_{\tau = 1}^{\infty} 
    g(j+\tau) u (-\tau) ,
\end{align}
where~$H_g$ is  called the 
 \emph{Hankel operator}
 corresponding to~$g$.
%
The variation-diminishing properties of such convolution operators is a central theme in the theory of total positivity~\cite{karlin_tp}, and has found applications in various fields including statistics and interpolation theory.

In matrix notation,
\begin{align}\label{eq:hankel_infi}
    \begin{bmatrix}
    y(0) \\
    y(1) \\
    y(2) \\
    \vdots \\
    \end{bmatrix}
    =
    \begin{bmatrix}
    g(1) & g(2) & g(3) & \hdots \\
    g(2) & g(3) & g(4) & \hdots \\
    g(3) & g(4) & g(5) & \hdots \\
    \vdots & \vdots & \vdots & \ddots
    \end{bmatrix}
    \begin{bmatrix}
    u(-1) \\
    u(-2) \\
    u(-3) \\
    \vdots \\
    \end{bmatrix}.
\end{align}
Thus,  the mapping from the input sequence to the output sequence is described by an infinite-dimensional Hankel matrix. Furthermore, for any~$u$ such that~$s^-(u),s^-(y)<\infty$,
there exists~$N>0$ such that the 
subsequences   $\{u(-N),u(-N+1),\dots u(-1)\}$ and~$\{y(0),y(1),\dots,  y(N-1)\}$  include  all the sign variations in~$u$ and~$y$, respectively, and the mapping between these subsequences is via a square  finite-dimensional Hankel matrix. 
Thus, studying an   SVDP  in this context can be done using the tools
described  in Section~\ref{sec:posi} (and some limit arguments like Prop.~\ref{prop:limitschanges}), but there are two twists:   (1)~we only need to consider Hankel matrices that are known to have   special total positivity  properties~\cite{Fallat_2017}; 
and   (2)~the SVDP can be related to 
the impulse response and/or  
transfer function of the system.

The Hankel operator~$H_g$ can also be expressed as
\[
            H_g u = O(A,c) (C(A,b) u),
\]
with
\begin{align*}
    C(A,b)u:&=\sum_{k=-\infty}^{-1} A^{-(k+1)} b u(k),\\
    (O(A,b) z)(j) :&=c^TA^j  z.
\end{align*}
Note that~$C(A,b)u$ is just~$x(0)$  i.e.,  the state of~\eqref{eq:dtisos} at time~$0$.

Any contiguous~$q\times q$ submatrix of the infinite Hankel matrix in~\eqref{eq:hankel_infi} is  of  the form 
\begin{align*}
H_{g} (p,q) := 
    \begin{bmatrix}
    g(p)     & g(p+1) & \hdots & g(p+q-1) \\
    g(p+1)   & g(p+2) & \hdots & g(p+q)   \\
    \vdots   & \vdots & \ddots & \vdots   \\
    g(p+q-1) & g(p+q) & \hdots & g(p+2q-2)
    \end{bmatrix},
\end{align*}
where $p,q\geq 1$.
The matrix~$H_{g} (p,q) $   can also
be associated with the state-space representation of the system.
Define the matrices 
$C^{p} (A,b)\in\R^{n\times p}$
and
$O^{p} (A,c)\in\R^{p\times n}$ by
\begin{align*}
    C^{p} (A,b) := 
    \begin{bmatrix}
    b & Ab & \hdots & A^{p-1} b
    \end{bmatrix},
    \quad 
     O^{p} (A,c) := 
    \begin{bmatrix}
    c^{T} \\
    c^{T} A \\
    \vdots \\
    c^{T} A^{p-1}
    \end{bmatrix}.
\end{align*}
%
%
Then a calculation shows that
\begin{align*}
    H_{g} (p,q) = O^{q} (A,c) A^{p-1} C^{q} (A,b).
\end{align*}
%

%
%
%
%

The SVDP property for the system is defined as follows. 
\begin{Definition}\label{def:OVDk}\cite{gruss2}
%
Fix~$k\geq0$. The system is called \emph{Hankel $k$-positive}  if
  for any 
sequence~$u$ with~$s^{-} (u) \le k-1$, we have 
\begin{enumerate}
    \item $s^{-} (y) \le s^{-} (u)$, and
    \item if $s^{-} (y) = s^{-} (u)$ then the signs of the first non-zero element in $u$ and in $y$ are the same.\label{item:CVDKITEM2}
\end{enumerate}
%
\end{Definition}

%
%
Note that here the requirement is   for any sequence~$u$ with~$s^{-} (u) \le k-1$,
and not only for sequences~$u$ such that~$s^{-} (u) =  k-1$. It should be clear from the discussion above and Remark~\ref{rem:thnr} that this requirement is closely related to~$TN_k$ of the Hankel matrix.

Recall that in general total positivity properties are not preserved under 
matrix sum (see Example~\ref{exa:not_sum}).
However, if~$A,B$ are both Hankel matrices and $TN_r$ then~$A+B$ is Hankel and~$TN_r$~\cite{Fallat_2017}. 
This can be used to prove the following. 
\begin{Proposition}\cite{gruss2}
If~$H_{g_1}$ and~$H_{g_2}$ are Hankel $k$-positive then so is~$H_{g_1}+H_{g_2}$.
\end{Proposition}
In other words, the parallel interconnection of systems preserves Hankel~$k$-positivity.

We now consider several  important  special cases of Definition~\ref{def:OVDk}. 
\subsection{Hankel $ 1$-positive  systems}
In a  Hankel $ 1$-positive   system,  
the output corresponding to   any sequence~$u$
with~$s^-(u)=0$ satisfies 
\[
s^-(y)=0,
\]
and the signs of the first non-zero element in~$u$ and in~$ y $ agree. 
Assume that at least one entry of~$u$ is positive. Then all the entries of~$u$ are non-negative, and all the entries of~$y $ are non-negative. Such a system is called~\emph{externally positive}, where the term externally refers to the fact  that 
the I/O mapping is positive,
but there are no requirements on the sign variations 
evolution in  the state vector~$x$. 

It is straightforward to characterize  Hankel $ 1$-positive  systems
in terms of the impulse response~$g$.
Let~$\delta(k)$ denote the discrete time impulse i.e.~$\delta(0)=1$, and~$\delta(k)=0$ for any~$k\not = 0$. 
Since~$s^-(\delta)=0$, a necessary condition for
Hankel $ 1$-positivity  
is that~$g(k)\geq 0$ for all~$k$. 
Eq.~\eqref{eq:conv} implies that  this condition is also sufficient.  

Given~$b,c\in\{0,1\}^n$ and~$A\in\{0,1\}^{n\times n}$, consider the sequence of non-negative numbers~$\gamma_k:=c^T A^ k b $, $k=0,1,\dots$. 
 Note that every~$\gamma_k$  satisfies either~$\gamma_k=0 $ or~$\gamma_k\geq 1$. 
It is known that determining if this sequence includes  a zero  is NP-hard in the dimension~$n$~\cite{BLONDEL200291}. This implies that given the state-space representation~\eqref{eq:dtisos}, it is NP-hard to determine if the system is Hankel $ 1$-positive.
 
  The next example illustrates the difficulty in determining  Hankel $ 1$-positivity  from
  the transfer function. 
\begin{Example}
%
Consider the \emph{first order lag}
\begin{align}\label{eq:fol}
    G(z) = \frac{r}{z-p}.
\end{align}
If~$p=0$ then~$g(k)= r\delta(k-1)$. If~$p\not =0$ then~$g(k)=r  p ^{k-1}\mathbbm{1}(k-1) $. Thus,~\eqref{eq:fol} is  Hankel $ 1$-positive   iff~$p,r\geq0$ i.e., iff~\eqref{eq:fol}  is a positive first-order lag. 

As another example, consider the  transfer function:
\begin{align*}
    G(z) = \frac{r_1}{z-p_1}
    +
    \frac{r_2}{z-p_2} , 
\end{align*}
where $p_1$ is a dominant pole, i.e.~$|p_1| > |p_2|$.
Since  
\begin{align*}
    g(k) = 
    ( r_1 p_{1}^{k-1} + r_2 p_{2}^{k-1} )
    \mathbbm{1}(k-1) , 
\end{align*}
  we have~$g(k)\approx
    r_1 p_{1}^{k-1}
$    for large values of~$k$. 
Thus, a necessary condition 
for Hankel $1$-positivity  is that~$p_1>0,r_1 \ge 0$, but this  is not a sufficient condition.
%
%
 
\end{Example}

\subsection{Hankel $ 2$-positivity    and unimodality}
The first-order difference of a sequence~$s$ is~$(\Delta^{(1)} s )(k): = s(k+1) - s(k)$, i.e. the ``discrete-time derivative'' of~$s$, and the~$j$th-order difference for~$j>1$ is defined inductively by
\[ (\Delta^{(j)} s )(k): =(\Delta^{(j-1)} s ) (k+1) - (\Delta^{(j-1)} s ) (k)  .\]

 The linearity of~\eqref{eq:dtisos} allows to 
 transform an   SVDP between~$u$ and~$y$ to an SVDP between~$\Delta^{(j)} u $
 and~$\Delta^{(j)} y $. 
For example, 
\begin{align}\label{eq:deltau}
    ( \Delta^{(1)} y )(k) &= y(k+1) - y(k) \nonumber\\
    &= ( g \circledast  u)(k+1) - ( g \circledast u )(k) \nonumber\\
    &= g  \circledast  ( u(k+1) - u(k) ) \nonumber\\
    &= g \circledast  (\Delta^{(1)} u) (k),
\end{align}
where~$\circledast$   denotes the convolution operator.

A 
sequence
$u $ 
is called \emph{unimodal} if 
$s^-(\Delta^{(1)} u )\leq 1$ (see, e.g.~\cite{unimodal}).
For example,~$\delta$ is unimodal, as 
\[
(\Delta^{(1)} \delta )(k)=\begin{cases}
        1, &k=-1,        \\
        -1, &k=0,        \\
        0,& \text{otherwise}. 
\end{cases}
\]
Eq.~\eqref{eq:deltau} implies that if~\eqref{eq:dtisos} is Hankel $2$-positive 
then it maps any unimodal input to a unimodal output. In particular, since~$\delta$ 
is unimodal, $g$ must be unimodal. 
 
 \subsection{Hankel $\infty$-positive systems}
Definition~\ref{def:OVDk} implies that for 
a Hankel $ \infty$-positive system
 every input/output pair satisfies
 $s^{-} (y) \le s^{-} (u)$,
 and if~$s^{-} (y) = s^{-} (u)$  then the signs of the first non-zero element in $u$ and in $y$ are the same. Such systems are also called Hankel totally positive or relaxation systems~\cite{relax_systems}, and have 
 special control-theoretic properties~\cite{rantzer_relax}. 
  These systems also  admit a simple characterization
 in  terms of their transfer function, namely, a system is Hankel $ \infty$-positive
 iff its  transfer function has the form 
 \[
 G(z)  = \sum_{i=1}^n \frac{r_i}{z-p_i}
 \]
with~$r_i, p_i \geq 0$ i.e., it is the parallel interconnection  of positive lags.

An important tool in the analysis of 
Hankel $k$-positivity  is the   Hankel     $k$-compound system. 
\begin{Definition}\cite{grussler2022variation}
Given an impulse response~$g$ and an integer~$k\geq 1$, the 
  \emph{Hankel $k$-compound system} is 
the system with impulse response
\[
g^{(k)}(j) =\det( H_g(j,k )  ),\quad  j=1,2,\dots 
\]
\end{Definition}
The transfer function associated with~$g^{(k)}$ is denoted by~$G^{(k)}$. 
   
\begin{Example}
The  Hankel $1$-compound system has impulse response
\[
g^{(1)}(j) =\det( H_g(j,1 )  )  =g(j)
\]
i.e., it is just the original system, and~$G^{(1)}(z)=G(z)$. 
The  Hankel $2$-compound system has impulse response
\begin{align*}
    g^{(2)}(j) &=\det( H_g(j,2 )  )  \\
    &=\det( \begin{bmatrix}
    g(j)&g(j+1) \\ g(j+1)& g(j+2)
    \end{bmatrix})\\
    &= g(j)g(j+2)-(g(j+1))^2. 
\end{align*}

\end{Example}

Recall from Prop.~\ref{prop:Fallat}
that in general verifying total positivity  of a matrix requires checking only contiguous and initial minors. In general, these results do not extend to checking total non-negativity. For example,
all the contiguous minors  of~$A=\begin{bmatrix}
0&1&0\\0&0&1\\1&0 &0
\end{bmatrix}$  are  non-negative,
yet~$A$  is not $TN_3$ as it admits non-contiguous minors that are negative. 
For the Hankel case,  things are simpler.
\begin{Proposition}\cite{gruss2}
        Let~$A\in\R^{n\times n}$ be Hankel. If all the contiguous minors of~$A$ up to order~$r$ are non-negative then~$A$ is~$TN_r$. 
\end{Proposition}

Using this, it is possible to relate Hankel~$k$-positivity  
to external positivity  of the Hankel $k$-compound system. 
\begin{Proposition}\cite{gruss2}
%
%
A DT SISO    LTI system with transfer function~$G$    
 is Hankel $k$-positive iff 
$G^{(k)}$ is externally positive.
%
\end{Proposition}
 
  In this respect,  the construction of~$G^{(k)}$ 
  is another demonstration of the general principal described in Section~\ref{sec:k_generalizations}.
%
\section{Discussion}
%

\begin{figure}
 \begin{center}
  \includegraphics[scale=0.8]{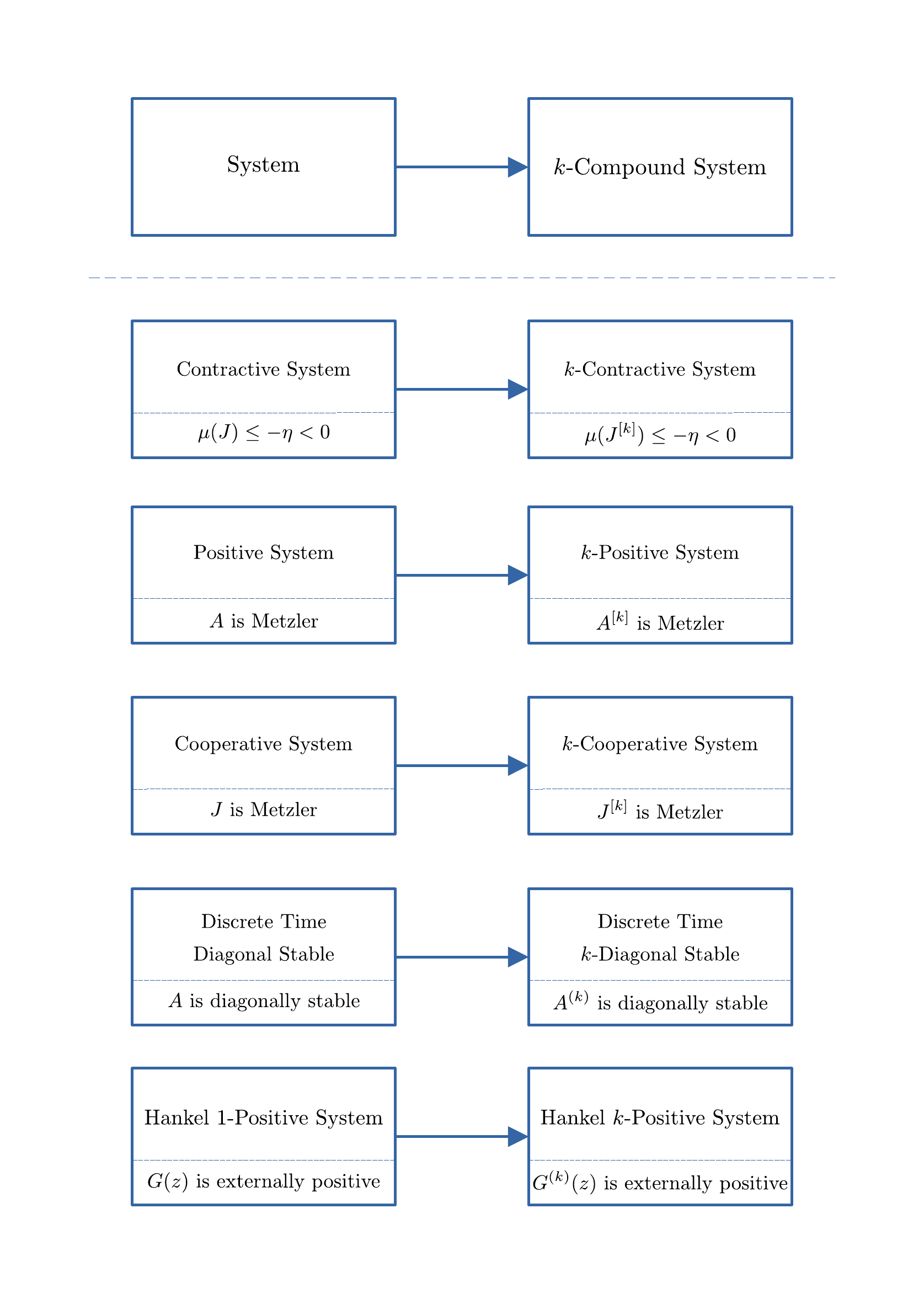}
	\caption{Generalizing various types of dynamical systems using  the $k$-compound system.}
	\label{fig:lift}
\end{center}
\end{figure}

We introduced a  principle that allows to generalize, in a non-trivial way,     various classes of linear and non-linear dynamical systems.
If the LTV system~\eqref{Eq:ltv}
  satisfies a specific \emph{property}, we  
   say that the~LTV satisfies~\emph{$k$-property} if the associated  $k$-compound system 
satisfies this property. This has been  used to generalize important classes of systems including contracting, positive, cooperative, and diagonally stable systems      into~$k$-contracting, $k$-positive, $k$-cooperative, and~$k$-diagonally stable systems (see Fig.~\ref{fig:lift}). 
This approach makes sense 
  for at least two reasons. First, for~$k=1$ the generalization reduces to the original class of systems,
  as 
  the~$1$-compound system is just  the original system. Second, the $k$-compound system  
  is based on~$k$-compound matrices,  and these matrices   play an important role in
  describing  the evolution of~$k$-parallelotopes and the evolution of sign changes along the dynamics. 
  
  In some cases, these generalizations lead to a natural hierarchy. For example, systems that are $k$-contracting   w.r.t. to   the~$L_1,L_2$ or~$L_\infty$ norms are  also~$\ell$-contracting
  w.r.t. to the same  norm for any integer~$\ell\geq k$. 
  
We believe that the ideas described in this tutorial paper may lead to   interesting research directions.
We now describe several possible directions.

%

%

\subsection{Linear Systems with Inputs}
In the context of systems and control theory, it is important  to consider systems with inputs and outputs. We   reviewed  one possible approach for analyzing an SVDP for such systems in Section~\ref{sec:hankel}. 
A possible alternative is to use 
  a  formula for analyzing the evolution of compounds in  systems with an input  from~\cite{MULD_SYSTEMS_WITH_INPUT}. Consider the set of time-varying 
  non-autonomous ODEs: 
\be\label{eq:nonauto}
\dot x^i(t)=A(t)x^i(t)+f^i(t),\quad i=1,\dots,k.
\ee
Then it is shown in~\cite{MULD_SYSTEMS_WITH_INPUT} that:
\be\label{eq:muldo}
\frac{d}{dt} \begin{bmatrix}
x^1&\dots&x^k
\end{bmatrix}^{(k)} =A^{[k]} \begin{bmatrix}
x^1&\dots&x^k
\end{bmatrix}^{(k)} +\sum_{i=1}^k 
\begin{bmatrix}
x^1&\dots&x^{i-1}& f^i&x^{i+1}&\dots &x^k
\end{bmatrix}^{(k)}.
\ee
Note that in the autonomous case, i.e. when~$f^i\equiv 0$  for all~$i $, this reduces to~\eqref{eq:expat}.

Many of the notions described in this paper  may   perhaps be extended to such  systems. 
For example, Angeli and Sontag  extended cooperative systems (more generally, monotone systems) to systems with inputs~\cite{mcs_angeli_2003}. It may be of interest to extend the notion of~$k$-positive systems and~$k$-cooperative systems to the case of systems with inputs.

As another example, a fundamental notion in systems and control theory   is input-to-state stability~(ISS)~\cite{ISS2020}. It may be of interest to consider a~$k$  generalization of this notion. To explain the basic idea, consider the LTI
\be\label{eq:lti2}
\dot x=Ax+1_2u,
\ee
with~$A=\begin{bmatrix}
-1&3\\3&-1
\end{bmatrix}$ and~$1_2=\begin{bmatrix} 1& 1 \end{bmatrix}^T$. Since~$A$ is not Hurwitz, this system is not ISS. Fix  two initial conditions~$a,b\in\R^2$, and let
$\Phi(t):=\begin{bmatrix}
x(t,a) & x(t,b)
\end{bmatrix}$. Applying~\eqref{eq:muldo} with~$k=2$ gives 
\begin{align*}
\frac{d}{dt} \det(\Phi)&= -2 \det(\Phi)+\begin{bmatrix}
x(t,a) & 1_2 u 
\end{bmatrix}^{(2)}+\begin{bmatrix}
1_2 u & x(t,b)
\end{bmatrix}^{(2)}\\
&= -2 \det(\Phi)+(x_1(t,a)-x_2(t,a)+x_2(t,b)-x_1(t,b))u,
\end{align*}
and substituting the solution of~\eqref{eq:lti2} gives
\begin{align*}
\frac{d}{dt} \det(\Phi)&= -2 \det(\Phi)+ (a_1-a_2+b_2-b_1)\exp(-4t)u.
\end{align*}
This suggests that $2$-dimensional volumes are in some sense input to state stable~(ISS) (see~\cite{ISS2020}), but with a gain function that depends on the initial conditions. An interesting  research  topic is formulating these notions rigorously, and studying their implications for interconnected systems.

 \subsection{$k$-Contraction on Networks} 
 There is considerable  interest in dynamical system over networks motivated in part  by applications in multi-agent systems, the electric grid, neuroscience, epidemiology,  
 and more~\cite{eger10,epidemics_on_nets,bullo_networks_book}. 
 
 Wang and Slotine~\cite{wang_slotine_2005} suggested an approach  for the 
 stability analysis and control  synthesis of networked 
 systems   using 
 contraction theory. It may be of interest to extend these results  using $k$-contraction and, more generally, to formulate a suitable notion of the $k$-compound system  of such networks and their associated connection graphs.  
 
 \subsection{Kernels Satisfying an SVDP and Applications to PDEs} 
We focused here on matrices satisfying an SVDP and its implications to dynamical systems described by ODEs. There is, however, a rich theory on kernels that satisfy an SVDP~\cite{karlin_tp}. These have found some applications in the analysis of  PDEs (see, e.g.,~\cite{fusco_PDE,smith_1990} and the references therein), but we believe that this field is still largely unexplored.

\section*{Acknowledgements}
It is written  in the Mishna  that
one who learns from his fellow a   single letter  must treat him with great honor.
The third  author would like  to thank  E. D. Sontag  for teaching him a whole language.
\section{Appendix: Proof of Thm.~\ref{thm:CB} }
We begin with an auxiliary result.
\begin{Lemma}\label{lem:calcdet}
Let~$A\in\mathbb{C}^{n\times m}$,  $B\in\mathbb{C}^{m\times n}$. 
Then
\[
\det(AB)=\begin{cases}
\sum_{\alpha\in Q(n,m)} A(\{1,\dots,n\} | \alpha)B(\alpha|\{1,\dots,n\}) ,& \text{ if }  n\leq m,\\
0 ,& \text{ if } n>m.
\end{cases}
\]
%
\end{Lemma}
%
\begin{IEEEproof}
%
Let~$D:=\begin{bmatrix}
 0&A\\B &I_m
\end{bmatrix} \in \C^{ (n+m)\times (n+m) }$.
Since 
\[
\begin{bmatrix}
I_n&-A\\0 &I_m
\end{bmatrix}
D
=\begin{bmatrix}
 -AB &0 \\B &I_m
\end{bmatrix},
\] we have 
\be\label{eq:eqdet}
\det (D) =\det( -AB)= (-1)^n\det(AB). 
\ee
 
The Laplace expansion of~$D$ using minors  that include the first~$n$ rows   gives:
\begin{align}\label{eq:explap}
    \det(D)=\sum_{\alpha \in Q(n,n+m)} D(\{1,\dots,n\} | \alpha )  \cof(D(\{1,\dots,n\} | \alpha)),
\end{align}
where~$\cof(D(\alpha|\beta)):=(-1)^{\alpha_1+\dots+\alpha_n+\beta_1+\dots+\beta_n} D((\{1,\dots,n+m\}\setminus  \alpha)|(\{1,\dots,n+m\}\setminus \beta ))$ is the \emph{cofactor} of~$D$ corresponding to~$(\alpha|\beta)$.
%
%
%
If~$m<n$ then
then any submatrix in the form~$D[\{1,\dots,n\} | \alpha ]$ includes a column of zeros, so~$D(\{1,\dots,n\} | \alpha ) =0$ and~$\det(D)=0$. Combining this with~\eqref{eq:eqdet} yields~$\det(AB)=0$. 

Now suppose that~$m\geq n$.
If~$\alpha$ includes any of the indexes~$1,\dots,n$ then~$D[\{1,\dots,n\} | \alpha]  $   includes a column of zeros. Using this and the fact that rows~$n+1,\dots,n+m$ of~$D$ includes rows of~$B$ and of~$I_m$ in~\eqref{eq:explap} gives
\[
\det(D)=(-1)^n  \sum_{\beta \in Q(n,m)}  
 A(\{1,\dots,n\} | \beta )  B (\beta |  \{1,\dots,n\} ),
\]
and combining this with~\eqref{eq:eqdet} completes the proof of Lemma~\ref{lem:calcdet}.
%
\end{IEEEproof}

We can now prove
Thm.~\ref{thm:CB}.
Let~$P\in\mathbb{C}^{n\times m}$,  $Q\in\mathbb{C}^{m\times s}$, and fix~$k\in \{1,\dots,\min \{n,m,s\} \}$. Let~$C:=PQ\in\C^{ n\times  s  }$. Fix~$\alpha \in Q(k,n)$ and~$\beta \in Q(k,s)$. Since~$c_{ij}=\sum_{\ell=1}^m p_{i\ell} q_{ \ell j}$, 
\[
C[\alpha|\beta]= P[\alpha|\{1,\dots,m\}] Q[\{1,\dots,m\} | \beta] .
\]
Let
\begin{align} \label{eq:abdef}
A&:=P[\alpha|\{1,\dots,m\}]\in 
\C^{  k\times m },\\
B&:=Q[\{1,\dots,m\} | \beta]\in
\C^{  m\times k },\nonumber
\end{align}
so that~$C[\alpha|\beta]=AB$.
Applying Lemma~\ref{lem:calcdet} gives
\begin{align*}
    C(\alpha|\beta)&= \det(AB)\\
&=\sum_{\gamma\in Q(k,m)} A(\{1,\dots,k\} | \gamma)B(\gamma|\{1,\dots,k\}),
\end{align*}
and using~\eqref{eq:abdef} gives
\begin{align*}
    C(\alpha|\beta)
&=\sum_{\gamma\in Q(k,m)} P(\alpha | \gamma)Q(\gamma|\beta).
\end{align*}
This completes the proof of Thm.~\ref{thm:CB}.
\bibliographystyle{IEEEtranS}
\bibliography{refs,compound, compound_bibtex,VDP}

 \begin{IEEEbiography}[{\includegraphics[width=1in,height=1.5in,clip,keepaspectratio]{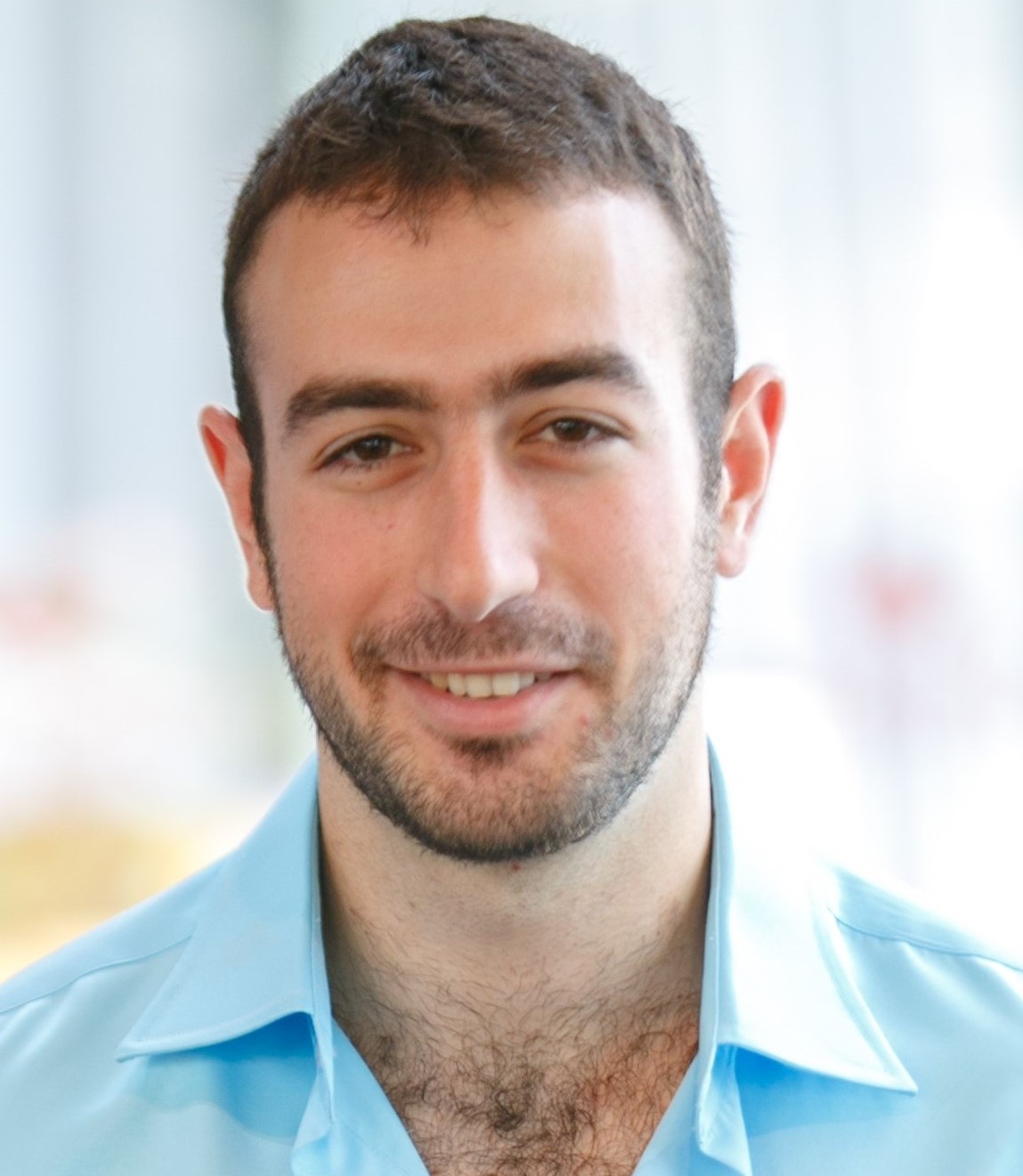}}]{Eyal Bar-Shalom} received the BSc degree (cum laude) and MSc degree (cum laude) in Elec. Eng. from Tel Aviv University, in 2012 and 2019, respectively.
He is currently a PhD student at the Dept. of Elec. Eng. at Tel Aviv University. 
\end{IEEEbiography}
  \begin{IEEEbiography}[{\includegraphics[width=1.1in,height=2in,clip,keepaspectratio]{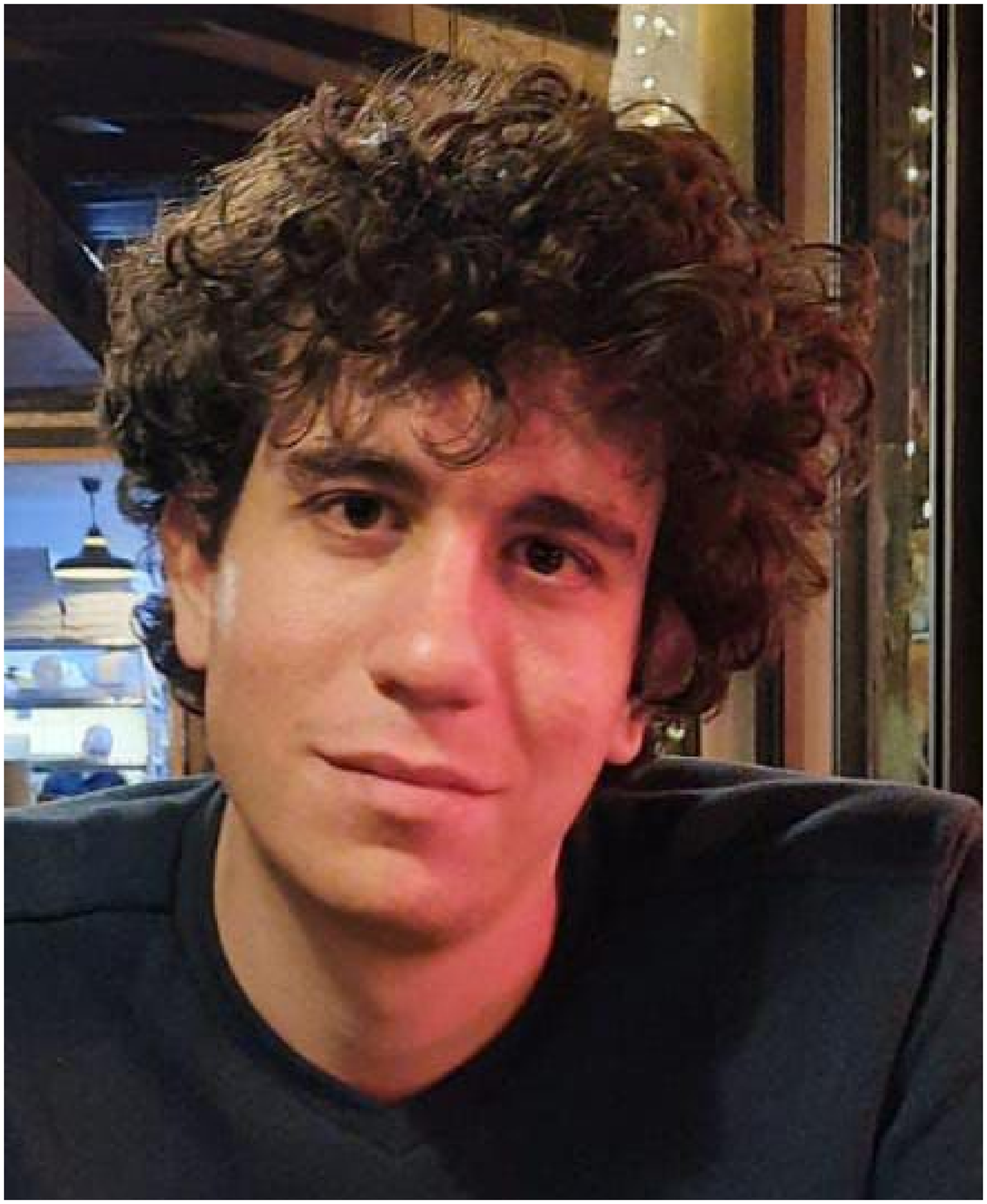}}\\]{Omri Dalin} received the  BSc degree in Mech. Eng. from Tel Aviv University, in 2020.
He is currently an MSc student at the Dept. of Mech. Eng. at Tel Aviv University. 
\end{IEEEbiography}
  \begin{IEEEbiography}
    [{\includegraphics[width=1in,height=1.25in,clip,keepaspectratio]{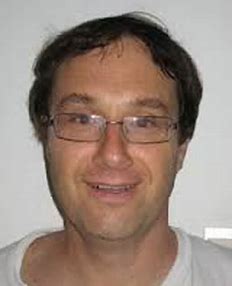}}]{Michael Margaliot}
 received the BSc (cum laude) and MSc degrees in
 Elec. Eng. from the Technion-Israel Institute of Technology-in
 1992 and 1995, respectively, and the PhD degree (summa cum laude) from Tel
 Aviv University in~1999. He was a post-doctoral fellow in the Dept. of
 Theoretical Math. at the Weizmann Institute of Science. In 2000, he
 joined the Dept. of Elec. Eng.-Systems, Tel Aviv University,
 where he is currently a Professor. His  research
 interests include the stability analysis of differential inclusions and
 switched systems, optimal control theory, computation with
 words, Boolean control networks, contraction theory, $k$-positive systems, and systems biology.
 He is co-author of \emph{New Approaches to Fuzzy Modeling and Control: Design and
 Analysis}, World Scientific,~2000 and of \emph{Knowledge-Based Neurocomputing}, Springer,~2009. 
 He  served
as  an Associate Editor of~\emph{IEEE Transactions on Automatic Control} during 2015-2017.
\end{IEEEbiography}

\end{document}